\documentclass[12pt]{article}
\usepackage[T1]{fontenc}
\usepackage[latin1]{inputenc}
\usepackage{amsmath,amsfonts,amssymb,amscd,mathrsfs,latexsym,epsfig,accents}
\usepackage{amsthm}
\usepackage{eufrak}
 \ifx\pdfoutput\@undefined\usepackage[usenames,dvips]{color}
 \else\usepackage[usenames,dvipsnames]{color}
 \IfFileExists{pdfcolmk.sty}{\usepackage{pdfcolmk}}{}
 \fi
\definecolor{anand}{rgb}{0,0.5,0}
\definecolor{error}{rgb}{0.8,0,0}
\definecolor{clifton}{rgb}{0,0,0.8}
\definecolor{k}{rgb}{.5,0,.5}

\newcommand{\C}{{\EuFrak C}}

\newcommand{\bi}{\begin{itemize}}
\newcommand{\ei}{\end{itemize}}

\def\Ind#1#2{#1\setbox0=\hbox{$#1x$}\kern\wd0\hbox to 0pt{\hss$#1\mid$\hss}
\lower.9\ht0\hbox to 0pt{\hss$#1\smile$\hss}\kern\wd0}
\def\ind{\mathop{\mathpalette\Ind{}}}
\def\Notind#1#2{#1\setbox0=\hbox{$#1x$}\kern\wd0\hbox to 0pt{\mathchardef
\nn=12854\hss$#1\nn$\kern1.4\wd0\hss}\hbox to
0pt{\hss$#1\mid$\hss}\lower.9\ht0 \hbox to
0pt{\hss$#1\smile$\hss}\kern\wd0}

\def\thind{\mathop{\mathpalette\Ind{}}^{\text{\th}} }

\def\nthind{\mathop{\mathpalette\Notind{}}^{\text{\th}} }

\def\starind{\mathop{\mathpalette\Ind{}}^{*} }

\def\nstarind{\mathop{\mathpalette\Notind{}}^{*} }

\newcommand {\thorn} {\text{\th}}

\newcommand{\acl}{\operatorname{acl}}
\newcommand{\tp}{\operatorname{tp}}
\newcommand{\dcl}{\operatorname{dcl}}

\def\uth{\text{U}^{\text{\th}} }

\makeatletter
\def\barx{\accentset{{\cc@style\underline{\mskip10mu}}}}

\def\barxv{\accentset{{\cc@style\underline{\mskip15mu}}}}

\def\barxx{\accentset{{\cc@style\underline{\mskip20mu}}}}

\def\barxxv{\accentset{{\cc@style\underline{\mskip25mu}}}}

\def\barxxx{\accentset{{\cc@style\underline{\mskip30mu}}}}

\def\barxxxv{\accentset{{\cc@style\underline{\mskip35mu}}}}
\makeatother

\newtheorem{theorem}{Theorem}[section]
\newtheorem{main theorem}{Theorem}
\newtheorem{proven main theorem}{Theorem}
\newtheorem{lemma}[theorem]{Lemma}
\newtheorem{fact}[theorem]{Fact}
\newtheorem{claim}{Claim}
\newtheorem{corollary}[theorem]{Corollary}
\newtheorem{proposition}[theorem]{Proposition}
\newtheorem{definition}[theorem]{Definition}
\newtheorem{remark}[theorem]{Remark}
\newtheorem{conjecture}[theorem]{Conjecture}
\newtheorem{question}[theorem]{Question}

\newtheorem{example}[theorem]{Example}

\title{Superrosy dependent groups having finitely satisfiable generics}
\author{Clifton Ealy, Krzysztof Krupi\'nski\footnote{Research supported by the Polish Government grant: N201 032 32/2231},  Anand Pillay\footnote{Research supported by a Marie Curie Chair: 024052}}
\date{}

\begin{document}
\maketitle
\begin{abstract}
We develop a basic theory of rosy groups and we study groups of small $\uth$-rank satisfying NIP and having finitely satisfiable generics: $\uth$-rank 1 implies that the group is abelian-by-finite, $\uth$-rank 2 implies that the group is solvable-by-finite, $\uth$-rank 2, and not being nilpotent-by-finite implies the existence of an interpretable algebraically closed field.
\end{abstract}

\section[\mbox{}]{Introduction}

\subsection*{Motivation}
If one hopes to apply geometric stability theory methods to study groups which are not necessarily stable, the weakest possible assumption that seems to be necessary, is rosiness, i.e. the assumption that we have an independence relation satisfying a minimal list of nice properties necessary to develop forking calculus. 

A general goal is to apply techniques from stable groups to the much wider class of rosy groups. During the last ten years, significant progress in the studies of groups in simple theories (which are always rosy) has been made. In this paper, following \cite{NIP paper}, we concentrate on another generalization of stable groups, namely, we will study rosy groups satisfying NIP and having finitely satisfiable generics (definitions to follow). In particular, our results generalize the appropriate theorems about stable groups and definably compact groups definable in o-minimal expansions of real closed fields. 

Another motivation is the fact that in the same way as groups of finite Morley rank generalize algebraic groups over algebraically closed fields, the groups that we consider (i.e. superrosy groups satisfying NIP and having finitely satisfiable generics) are a common generalization of algebraic groups over algebraically closed fields and compact Lie groups.

There is a long history of structural theorems about groups in model theory.  For instance, one has the following \cite{superstable fields}:

\begin{theorem}\label{superstable structure} Let $G$ be a superstable group. Then\\
(a) if it is of $U$-rank 1, it is abelian-by-finite, and\\
(b) if it is of  $U$-rank 2, it is solvable-by-finite.
\end{theorem}

\medskip
\noindent In comparison, here are the corresponding theorems about groups in o-minimal structures \cite{one dim o-min group, two and three dim o-min groups}:


\begin{theorem} Let $G$ be a definably connected group definable in a o-minimal structure.  Then

(a) (Razenj) if it is of dimension 1, it is isomorphic to either $\bigoplus_{p\in P}\mathbb Z_{p^\infty}\oplus\bigoplus_\delta\mathbb Q$ or $\bigoplus_\delta\mathbb Q$, and, in particular, it is abelian.
 
(b) (Nesin, Pillay, Razenj) if it is of dimension 2, and non-abelian, it is $R^{+}\rtimes R^{*}_{>0}$, for some real closed field $R$, and, in particular, it is solvable.

\end{theorem}

\medskip

Our initial goal was to obtain results along these lines in the more general setting of superrosy groups with hereditary fsg and the NIP.  Indeed, we have proved the following theorems.

\begin{main theorem}\label{Main Theorem 1}
Let $G$ be superrosy with NIP. Then if $G$ is of  $\uth$-rank 1 and has fsg, it is abelian-by-finite.
\end{main theorem}

\begin{main theorem} \label{Main Theorem 2}
Let $G$ be superrosy with NIP. Then if $G$ is of   $\uth$-rank 2  and has hereditarily fsg,  it is solvable-by-finite.
\end{main theorem}

While in the case of groups, there are a rich class of non-stable examples of dependent, superrosy groups with finitely satisfiable generics, the same is not true of fields.  In fact, we show

\begin{main theorem} \label{Main Theorem 3}
Suppose that $K$ is a superrosy  field and $K^+$ has fsg. Then $K$ is algebraically closed.
\end{main theorem}

\noindent Finally, we show 
\begin{main theorem} \label{Main Theorem 4}
Assume that $G$ has  NIP, hereditarily fsg, $\uth(G)$=2 and $G$ is not nilpotent-by-finite. Then, after possibly passing to a definable subgroup of finite index and quotienting by its finite center, $G$ is (definably) the semidirect product of the additive and multiplicative groups of an algebraically closed field $F$  interpretable in $G$, and moreover $G = G^{00}$.
\end{main theorem}

\subsection*{Preliminaries}
Throughout, we may assume that we work in a  monster model (i.e. large $\kappa$-saturated model) $\mathfrak C$ of a theory $T$. $G$ will always denote a definable group in this model, and all collections of parameters, $A,B,C$ etc., are assumed to be of size less than $\kappa$.  (With the one exception, of course, being when we consider a global type, $\bf p$, over all of $\C$.)  We write $G$ multiplicatively, with identity $e$.

When we write ``definable'' or ``type definable'' without further qualification, we will always mean ``definable with parameters'' or ``type definable with parameters'' respectively.

\subsubsection*{Rosy Theories}
 A model, $M$, is said to be rosy if it admits a notion of independence which extends to $M^{eq}$.  More precisely, we have the following definition (See, for instance, \cite{properties and consequences of thorn-independence}, or \cite{forking and thorn forking} for an alternate treatment):

\begin{definition}\label{rosy}
T is \emph{rosy} if and only if there exists a ternary relation, $\starind$, on both real and imaginary subsets (we treat tuples as subsets) of models of T such that:

\begin{enumerate}
\item $\starind$ is automorphism invariant.
\item If $c\in$acl$(aB)\setminus$acl$(B)$, then $a\nstarind_B c$.
\item If $a\starind_B C$ and $B \cup C \subseteq D$, then there is some $\tilde a\models \tp(a/BC)$ such that $\tilde a\starind_B D$.
\item There is some $\lambda$ such that for all $a$, whenever one has $(B_i)_{i<\alpha}$ with $B_i\subset B_j$ for $i<j$ and $a\nstarind_{B_i} B_j$ for $i<j<\alpha$, then $\alpha<\lambda$.
\item If $B\subseteq C \subseteq D$, then $a\starind_B D$ if and only if $a\starind_B C$ and $a\starind_C D$.
\item $C\starind _A B $ if and only if $c\starind _A B$ for any
finite $c\subseteq C$.
\item $a\starind_C b$ if and only if $b\starind_C a$.
\end{enumerate}

\noindent We refer to a relation that satisfies (1) to (7) an \emph{independence relation}. If $a\nstarind_C b$, we say that $\tp(a/Cb)$ \emph{$*$-forks} over $C$. A model is \emph{rosy} if its theory is rosy. 
\end{definition}
   
\begin{remark}
Any theory with a ternary relation on  $\C^{eq}$ satisfying (1) to (4) together with the left to right direction of (5) is rosy, but these need not imply that the ternary relation is an independence relation.  We restrict our attention to $\starind$.
\end{remark}

\noindent Alternatively, in any rosy theory, one may define a particular, well-behaved, notion of independence, namely \th-forking, and give an equivalent definition of rosiness based on the behaviour of this notion of independence \cite{characterizing rosy theories}. 

\begin{definition}

A formula $\delta (x,a)$ \emph{strongly divides over $A$} if the formula is not almost over $A$ and $\{ \delta (x,a') \}_{a'\models
\tp(a/A)}$ is $k$-inconsistent for some $k\in \mathbb{N}$.

We say that $\delta (x,a)$ \emph{\th-divides over $A$} if we can
find some tuple $c$ such that $\delta (x,a)$ strongly divides over
$Ac$.

A formula \emph{\th-forks over $A$} if it implies a (finite)
disjunction of formulas which \th -divide over $A$.

We say that the type $p(x)$
\emph{\th-divides over $A$} if there is a formula in $p(x)$ which
\th-divides over $A$; \th-forking is similarly defined. We say
that $a$ is \emph{\th-independent from $b$ over $A$}, denoted
$a\thind_A b$, if $\tp\left( a/Ab\right) $ does not \th -fork over
$A$.
\end{definition}

\begin{fact}
 A theory is rosy if and only if \th-forking is a symmetric relation.
\end{fact}

Not only is \th-forking symmetric, in fact we have (e.g. \cite{properties and consequences of thorn-independence}):

\begin{fact}
In a rosy theory, \th-independence is an independence relation.
\end{fact}

Of all independence relations, \th-independence is the weakest in the following sense:

\begin{fact}
Let $\starind$ be an independence relation (or, in fact, any ternary relation satisfying (1)-(4) and the left to right direction of (5)), then $a\starind_C b \Rightarrow a\thind_C b$. 
\end{fact}

\begin{example}
1. Any simple theory is rosy.\\
2. Any o-minimal theory is rosy.\\
3. The theory of p-adically closed valued fields is not rosy.
\end{example}

The structures in which we are interested are not just rosy, but \emph{superrosy}.  We will define superrosy theories in terms of the $\uth$-rank \cite{properties and consequences of thorn-independence}.

\begin{definition}
We define the \emph{$U^*$-rank} inductively as follows. Let $p$ be a
complete type over some set $A$. Then,
\begin{enumerate}
\item $U^*(p)\geq 0$ if $p$ is consistent. 

\item For any
ordinal $\alpha$, $U^*(p)\geq \alpha+1$ if there is some tuple
$a$ and some type $q \in S(Aa)$ such that $q\supset
p$, $U^*\left( q \right) \geq \alpha$ and $q$
$*$-forks over $A$. 

\item For any $\lambda$ limit ordinal,
$U^*(p)\geq \lambda$ if and only if $U^*(p)\geq \sigma $ for all
$\sigma < \lambda$.
\end{enumerate}
\end{definition}

\noindent  An important property of $U^*$ is that it satisfies the Lascar inequalities:

\begin{proposition}\label{lascar}
    $U^*(a/bA)+U^*(b/A)\leq U^*(ab/A)\leq U^*(a/bA)\oplus U^*(b/A)$
\end{proposition}
\noindent
{\em Proof.} As Theorem 5.1.6 in \cite{simple theories}. \hfill $\blacksquare$

\begin{definition}
T is \emph{superrosy} if and only if $\uth(p)<\infty$ for
every type $p$.
\end{definition}

Clearly, a theory is superrosy if and only if there is some independence relation such that that $U^*(p)<\infty$  for every $p$ (recalling, of course, that we insist that an independence relation extends to  $\C^{eq}$).

 \begin{definition}
For a definable set $X:=\varphi(\C,a)$ we define $U^*(X):=sup\{ U^*(p): p \in S(a), \varphi(x,a) \in  p\}$.
\end{definition}

If $T$ is of finite $U^*$-rank, then for every $X$ as above, $U^*(X)=U^*(p)$ for some $p \in S(a)$ containing $\varphi(x,a)$. The same conclusion is also true when $X$ is a definable group, as we will show later.

\begin{example}
1. Any supersimple theory is superrosy, and $\uth$-rank equals $SU$-rank.\\
2. Any o-minimal theory is superrosy, and $\uth$-rank equals dimension.
\end{example}

Finally, the following lemma about rosy theories is quite useful: namely that although the definition of \th-dividing requires that one produce a $k$-inconsistent set, to get \th-forking one need only find an \emph{almost} $k$-inconsistent set.  That is one need only to find a uniform family of formulas so that the intersection of any $k$ members of the family is finite.  There seems to be no reference for this fact, so we give a proof below.

\begin{lemma}\label{finite intersections good enough}
Suppose that $\varphi(a,b)$, and suppose that for any $k$ distinct realizations, $b_1,\dots, b_k$, of the non-algebraic $\tp(b/C)$ the conjunction $\bigwedge_{i\leq k}\varphi(x,b_i)$ defines a finite set, but $a\notin \acl(C)$.  Then $a\nstarind_C b$.
%
\end{lemma}
\noindent
{\em Proof.}
Suppose that for any $k$ distinct realizations, $b_1,\dots, b_k$, of $\tp(b/C)$ the conjunction $\bigwedge_{i\leq k}\varphi(x,b_i)$ defines a finite set. By compactness, there is $m$ such that each such conjunction defines a set of size less than $m$. Thus there is a maximal $n$ (possibly $n=0$) such that there is $\vec a=(a_1,\dots, a_n)$ for which one may find infinitely many $b_i\models\tp(b/C)$ such that each $\varphi(\C,b_i)$ contains $\vec a$.  We may assume, after possibly moving $\vec{a}$ by a $C$-automorphism, that $\tp(b/C\vec a)$ is not algebraic and $\varphi(\C,b)$ contains $\vec a$.

By the maximality of $n$, we see that for every $a' \in \varphi(\C,b)\setminus \{a_1,\dots,a_n\}$, $b \in \acl_C(a', \vec{a})$. If $a \in \acl(bC)$, we are done. Otherwise, we can choose $a' \models \tp(a/bC)$ different from $a_1, \dots, a_n$ so that 
\begin{enumerate}
\item[$(*)$] $a' \starind_{bC} \vec{a}$.
\end{enumerate} 
Since $b \in \acl_C(a',\vec{a})$ and $b \notin \acl_C(\vec{a})$, we get $b \nstarind_{C\vec{a}}a'$. This together with $(*)$ implies $a' \nstarind_C b$. Hence $a \nstarind_C b$. \hfill $\blacksquare$

\subsubsection*{Dependence}
In addition to rosiness, the groups we will consider satisfy a second condition: \emph{dependence}  (often called the non independence property and denoted by NIP):

\begin{definition}
    $T$ is said to have the \emph{independence property} if there is $\varphi(x,y)$, an $\mathfrak M\models T$,  and an infinite set $A\subseteq \mathfrak M$ such that for any $B\subseteq A$ one can find a $c_B$ such that $\varphi(\mathfrak M, c_B)\cap A = B$.  Otherwise, $T$ is said to be \emph{dependent}.
  \end{definition}

The only consequence of dependence in groups that we will use is the following, from \cite{logical stability in group theory}.

\begin{fact}\label{intersection is intersection of n}
    If $G$ is a group defined in a dependent theory, for each $\varphi$ there is some $n$ such that the intersection of any finite family of $\varphi$-definable subgroups is an intersection of $n$ members of the family. 
  \end{fact}

\begin{example}
        1. Any stable theory is dependent. In fact, simple and dependent is the same as stable, and supersimple and dependent is the same as superstable.\\
        2. $p$-adically closed fields are dependent.\\
        3. Pseudo-algebraically closed but not algebraically closed fields are not dependent.
\end{example}

\subsubsection*{Finitely Satisfiable Generics}
The final condition satisfied by the groups we consider is that they (and  sometimes their definable subgroups) have \emph{finitely satisfiable generics}.

First let us give a precise definition of translations of formulas and types which we will use in this paper.

\begin{definition}
We assume that $G$ is a group definable in $\C$ by a formula $G(x)$. Let $g \in G$, $\varphi(x,y)$ be any formula, $\pi(x)$ any partial type (over a small set) and ${\mathbf p(x)}$ any global type containing $G(x)$. We define:\\
1. $\varphi^*(x,wy):= (\exists u) (\varphi(u,y) \wedge G(w) \wedge G(u) \wedge x=w\cdot u)$ and $g\varphi(x,y):= \varphi^*(x,gy)$,\\
2. $g\pi(x):=\{g\varphi(x): \varphi(x) \in \pi(x) \}$,\\
3. $g{\mathbf p(x)}$ is the unique global type implied by $\{g\varphi(x): \varphi(x) \in {\mathbf p(x)} \}.$
\end{definition}

It is obvious that $g\varphi(x)$ and $g\pi(x)$ define the sets $g\cdot \varphi(G)$ and $g\cdot \pi(G)$, respectively. If $p \in S(A)$ extends $G(x)$ and $g \in A$, then $gp$ implies the unique complete type over $A$, and in some situations $gp$ will denote this complete type.

\begin{definition}
    A formula $\varphi(x)$ (or the set $\varphi(G)$) is \emph{left generic} if there are $g_1, \dots, g_n \in G$ such that  $g_1\varphi(G) \cup \dots \cup g_n\varphi(G) =G$. We say that a type is left generic if every formula in it is left generic.
\end{definition}

\begin{definition}
    $G$ has \emph{finitely satisfiable generics} (or \emph{fsg}) if there is a global type $\mathbf p$ containing $G(x)$ and a model $M\prec \mathfrak C$, of cardinality less than the degree of saturation of $\mathfrak C$, such that for all $g$, $g \mathbf p$ is finitely satisfiable in $M$ (i.e. each formula in $g \mathbf p$ defines a set which intersects $M$).
    
    We say that $G$ has \emph{hereditarily fsg} if every definable subgroup of $G$ also has fsg.

    When we consider a definable subgroup $G$ of some non-saturated model (as in the example below), we say that $G$ has (hereditarily) fsg if the set defined by the same formula in a large saturated model of theory has (hereditarily) fsg.
\end{definition}

\begin{remark}
    It is not difficult to see that if $Y\subseteq G$ is generic, and $N$ is a normal definable subgroup, then the image of $Y$ under the quotient map is generic in $G/N$.  Thus, if $G$ has fsg, and $N$ is a normal definable subgroup, then $G/N$ has fsg.  
\end{remark}

\begin{example}
        1. Algebraic groups have hereditarily fsg, as in fact do all stable groups, or even stably dominated groups.\\
        2. Compact Lie groups (which are interpretable in o-minimal structures) have hereditarily fsg, by \cite{NIP paper}.\\
	3. The complex numbers, in the language of fields together with a predicate defining the algebraic closure of rational numbers, has fsg, but not hereditarily fsg.
\end{example}

The following fact (from \cite{NIP paper}) is the central result about groups with finitely satisfiable generics:

\begin{fact}\label{basic consequences of fsg}
    Suppose that $G$ has fsg as witnessed by $\mathbf p$.  Then\\
    1. A formula is left generic if and only if it is right generic  (so we will skip the words `left' and `right').\\
    2. $\mathbf p$ is generic.\\
    3. The family of nongeneric subsets of $G$ form an ideal, so every partial generic type can be extended to a global one.\\
    4. $G^{00}$ (i.e. the smallest type definable subgroup of bounded index in $G$) exists, is type definable over empty set, and is the stabilizer of every global generic type.
  \end{fact}  
  
We mention here a couple of additional facts to be used later:

\begin{remark}
    Note that for every $X \subseteq G$ if $g_1 X\cup \dots \cup g_n X$ covers $G$, then $X^{-1}g_1^{-1} \cup \dots \cup X^{-1}g_n^{-1}$ covers $G$ as well.  Thus $X$ generic implies $X^{-1}$ is generic.
\end{remark}

\begin{proposition}\label{unique generic}
    Assume $G$ has fsg. If $G^{00}$ is definable, it has a unique global generic type.   
\end{proposition}
\noindent
{\em Proof.}
We may assume $G=G^{00}$. By Fact \ref{basic consequences of fsg}, at least one global generic type exists. Denote it by $\mathbf p$. Take any generic formula $\varphi(x)$. It is enough to show that $\varphi(x) \in {\mathbf p}$. 

There is $g\in G$ such that $g\varphi(x) \in {\mathbf p}$. So $\varphi(x) \in g^{-1}{\mathbf p}$. But by Fact \ref{basic consequences of fsg}, $G=G^{00}$ stabilizes $\bf p$, i.e. $g^{-1}{\mathbf p}={\mathbf p}$. So $\varphi(x) \in {\mathbf p}$. \hfill $\blacksquare$

\section{Rosy groups}

The facts about rosy groups we will use are all quite straightforward, and proofs in general follow the proofs about groups in simple and stable theories.  However, there is no prior exposition of these facts in the case of rosy groups, so we will provide proofs, or in the case where the proof is identical to that in simple theories, a reference.

We recall the definition of local \th-ranks, which we will use briefly in \ref{ucc} below and then repeatedly during our proof of the existence of \th-generics.

Throughout this section, $G$ will denote a group definable by a formula $G(x)$ (over $\emptyset$) in a monster model $\C$ of a rosy theory  $T$. Such a group will be called a rosy group.

\begin{definition}
Given a formula $\psi(x)$, a finite set $ \Phi $ of formulas
with object variables $x$ and parameter variables $y$, a finite set of
formulas $\Theta$ in the variables $y,z$, and natural number $k>0$, we define the \emph{\th$_{\Phi,\Theta, k}$-rank} of $\psi$  inductively as follows:
\begin{enumerate}

\item $\thorn_{\Phi,\Theta,k}(\psi) \geq 0$ if $\psi $
is consistent.

\item  For $\lambda $ limit ordinal, $\thorn_{\Phi,\Theta,k}(\psi)\geq \lambda $ if and only if $\thorn_{\Phi,\Theta,k}(\psi)\geq \alpha $ for all $\alpha
<\lambda $.

\item $\thorn_{\Phi,\Theta,k}(\psi) \geq \alpha +1$ if
and only if there is a $ \varphi \in \Phi $, some $\theta(y,z) \in
\Theta$ and parameter $c$ such that

\begin{enumerate}
\item  $\thorn_{\Phi,\Theta,k}(\psi\land \varphi(x,a))\geq \alpha $ for infinitely many $a\models \theta(y,c) $, and

\item  $\left\{ \varphi \left( x,a\right) \right\} _{a \models \theta(y,c)}$ is
$k-$inconsistent.
\end{enumerate}
\end{enumerate}

\noindent Given a (partial) type $\pi(x)$ we define $\thorn_{\Phi,\Theta,k}(\pi(x))$ to be the minimum of $\thorn_{\Phi,\Theta,k}(\psi)$ for $\psi\in\pi(x)$.  When $\Phi$ and $\Theta$ each contain only one formula, we will write $\thorn_{\varphi,\theta,k}(\psi)$.

\end{definition}

We recall that a theory is rosy if and only if for each $\psi, \Phi, \Theta, k$, the local thorn rank $\thorn_{\Phi,\Theta, k}(\psi)$ is finite.  Given a partial type $\pi(x)$, and a $\thorn_{\Phi,\Theta, k}$-rank, one can always extend  $\pi(x)$ to a complete type of the same $\thorn_{\Phi,\Theta, k}$-rank  as $\pi(x)$.   Moreover, $a\thind_C b$ if and only if for each $\Phi, \Theta, k$, one has that $\thorn_{\Phi,\Theta, k}(\tp(a/bC))=\thorn_{\Phi,\Theta, k}(\tp(a/C))$.

First, and easiest, we have a collection of chain conditions.  Definitions, (from, e.g., \cite{stable groups}) are repeated here for convenience.

\begin{definition}
    A family of groups $\{H_i:i\in I\}$ is called \emph{uniformly definable} if there is a formula $\varphi$ and parameters $b_i$ such that $H_i$ is defined by $\varphi(x,b_i)$.
    
    A group satisfies the \emph{uniform chain condition}, or \emph{ucc}, if for any formula $\varphi$ there is an $m_\varphi<\omega$ such that each chain of $\varphi$-definable groups has length at most $m_\varphi$.
\end{definition}

\begin{proposition}\label{ucc}
    A rosy group has the ucc.
\end{proposition}
\noindent
{\em Proof.}
    Note that if there is a chain of length $\omega$ of groups $H_i:=\varphi(\mathfrak C,b_i)$, then by compactness, there is a descending chain of arbitrary length.  Let $\{H_i:i<\omega\cdot\omega\}$ be such that for $i>j$, $H_i<H_j$.  Then for $k<\omega$, $[H_{k\omega}:H_{(k+1)\omega}]\geq \aleph_0$.  Let $\widetilde \varphi(x,b,d)$ be the formula that says ``$x$ is in the coset of $\varphi(\mathfrak C,b)$ that corresponds to the element $d$ of $M^{eq}$.''  This shows that $G$ has infinite  $\thorn_{\widetilde \varphi, \theta, 2}$-rank for appropriate $\theta$. \hfill $\blacksquare$

\begin{proposition}
    A superrosy group has the $\omega$dcc, i.e. any descending chain of definable groups, each with infinite index in its predecessor, is finite. 
\end{proposition}
\noindent
{\em Proof.} Suppose such a descending chain $G=H_0>H_1>\dots$ exists. Let $A$ be such that each $H_i$ is $A$-definable. 
We can choose a sequence $(a_i: i \in \omega)$ such that 
\begin{enumerate}
\item $b_0 \notin \acl(A)$, $b_1 \notin \acl_A(b_0)$, $b_2 \notin \acl_A(b_0,b_1)$, $\dots$, where $b_i \in \C^{eq}$ is a name of the coset $a_{i}H_{i+1}$,
\item $a_0H_0 \supseteq a_1H_1 \supseteq a_2H_2 \supseteq \dots$.
\end{enumerate}
Now take any $a \in \bigcap_{i \in \omega} a_iH_i$. 
Then $b_i \in \acl_A(a) \setminus \acl_A(b_{<i})$. So $a \nthind_{Ab{<i}}b_i$ for $i \in \omega$, contradicting superrosiness. \hfill $\blacksquare$\\[3mm]
\noindent The proof also shows the following:

\begin{proposition}
    A superrosy group of finite $\uth$-rank also satisfies the $\omega$acc: any ascending chain of definable groups, each with infinite index in its  succesor, is finite.
\end{proposition}

\begin{definition}
    A group has the \emph{intersection chain condition}, or \emph{icc}, if for each $\varphi$ there is some $m_\varphi<\omega$ such that any chain of intersections of $\varphi$-definable groups has length at most $m_\varphi$.
\end{definition}

\begin{proposition}\label{icc}
    A rosy dependent group has the icc.
\end{proposition}
\noindent
{\em Proof.}
   Dependent means we may apply Fact \ref{intersection is intersection of n}, and rosy means we may apply the ucc.  Together, they clearly give the icc. \hfill $\blacksquare$\\[3mm]
\noindent We never use the full icc in any proof.  Rather we use the following:
\begin{corollary}\label{icc on centralizers}
    Rosy dependent groups have the icc on centralizers.
\end{corollary}

\noindent Of particular use is the centralizer connected component of $G$.

\begin{definition}
    The \emph{centralizer connected component} is the intersection of all centralizers of finite index.
\end{definition}

\begin{corollary}\label{centralizer connected component} 
The centralizer connected component of a rosy dependent group has finite index in this group and is $\emptyset$-definable.
\end{corollary}

\subsection*{\th-Generics}

We define now a different notion of generic for a given independence relation $\starind$, modeled after the notion of a generic type in simple theories.  We first introduce a definition of a $*$-generic type, and prove several facts about them assuming that they exist, and then show that they exist in the particular case where the independence relation is $\thind$.  Most proofs are simply obtained from the standard proof in simple theories by replacing $\ind$ by $\starind$.  In these cases, we give a reference rather than a proof.

It is important to note that we deal with two different notions, both regularly called ``generic''.  If we say simply ``generic'', we will always mean generic in the sense that finitely many translates cover $G$.  The notion of genericity that arises from an independence relation $\starind$, which we are about to introduce, will always be referred to as $*$-generic.

\begin{definition}
We say that a type, $p \in S(A)$,  extending $G(x)$ is \emph{left $*$-generic over A} if for all $a,b \in G$ with $a\models p$ and $a\starind_A b$, one has that $b\cdot a\starind Ab$.  We say that it is \emph{right $*$-generic over A} if, for $a,b$ as above, we have $a\cdot b\starind Ab$.  A type is \emph{$*$-generic} if it is both right and left $*$-generic.  
\end{definition}

\begin{lemma} \label{basic facts about *-generics}
\begin{enumerate}
    \item If $p$ is left (right) $*$-generic then $p$ does not $*$-fork over the empty set.
    \item  Let $a, b \in G$. If $\tp(a/A)$ is left $*$-generic and $b \in \acl(A)$, then $\tp(b \cdot a/A)$ is also left $*$-generic.

    \item Let $p(x)$ be a type containing $G(x)$, and let $q$ be a non-$*$-forking extension of $p$.  Then $p$ is left $*$-generic if and only if $q$ is left $*$-generic.
    \item If $p\in S(A)$ and $B\subseteq A$, and $p$ is $*$-generic, then so is $p|_B$.
    \item If $p \in S(A)$ is left $*$-generic, then $p^{-1}$ is as well, where we define $p^{-1}$ to be $\tp(a^{-1}/A)$ where $a$ is any realization of $p$.
    \item A type is left $*$-generic if and only if it is right $*$-generic.
    \item If $a\starind_A b$ and $\tp(a/A)$ is $*$-generic, then so is $\tp(b\cdot a/A)$.

\end{enumerate}
\end{lemma}
\noindent
{\em Proof.}
1. Let $b=e$.\\
2. Lemma 4.1.2.1 of \cite{simple theories}.\\
3. Lemma 4.1.2.2 and Lemma 4.1.2.3 of \cite{simple theories}.\\
4. Left (right) $*$-generic types do not $*$-fork over the empty set.\\
5. Lemma 4.1.2.4 of \cite{simple theories}.\\
6. Note that $p$ is left $*$-generic if and only if $p^{-1}$ is right $*$-generic.\\
7. Note that $\tp(a/Ab)$ is $*$-generic, and thus so is $\tp(b\cdot a/Ab)$. From this we see that $\tp(b\cdot a/A)$ is also $*$-generic. \hfill $\blacksquare$\\

Now we prove the existence of $\thorn$-generic types.  As in the corresponding proof in the case of simple theories, the existence of $\thorn$-generics will follow from the definition and examination of translation invariant local ranks. 

We introduce the translation invariant stratified \th-ranks as a specific kind of local \th-rank:

\begin{definition}
    The \emph{stratified $\thorn_{\Phi, \Theta, k}^G$-rank} of a partial type $\pi(x)$  is defined as $\thorn_{\Phi^*, \widetilde\Theta, k}(\pi(x) \cup \{ G(x)\})$, where $\Phi^*=\{\varphi^*(x,wy):\varphi(x,y)\in \Phi\}$ and  $\widetilde\Theta(wy,w'z)=\{\theta(y,z) \land w=w' \land G(w):\theta\in \Theta\}$. 
\end{definition}

While the stratified ranks are defined for formulas, it is clear that if two different formulas define the same set, they have the same rank, so we may speak of the stratified ranks of definable sets as well.

\begin{lemma}
    The stratified rank is translation invariant.
\end{lemma}
\noindent
{\em Proof.}
    By induction on rank:  It is clear that for every  $\psi(x)$, for every $\Phi, \Theta$, and $k$, and for every $g\in G$, we have $\thorn_{\Phi, \Theta, k}^G(\psi(x))\geq 0$ if and only if  $\thorn_{\Phi, \Theta, k}^G(g\psi(x))\geq 0$.
    
    Now suppose that $\widetilde\theta(wy,hc)\in \widetilde\Theta$ and $\varphi\in\Phi$ witness that the $\thorn_{\Phi, \Theta, k}^G$-rank of $\psi(x)$ is greater than  $n$. That is, suppose 
    
   \begin{enumerate}
   \item  $ \thorn_{\Phi, \Theta, k}^G(\psi(x) \land \varphi^*(x,hb))\geq n$ for infinitely many $hb \models \widetilde\theta(wy,hc)$, and

   \item  $\left\{ \varphi^* \left( x,hb\right) \right\} _{hb \models \widetilde\theta(wy,hc)}$ is $k-$inconsistent.
   \end{enumerate}
   
    \noindent Then $\widetilde\theta(wy,(g\cdot h)c)$ and $\varphi \in \Phi$ witness that the $\thorn_{\Phi, \Theta, k}^G$-rank of  $g\psi(x)$ is greater than $n$. \hfill $\blacksquare$

\begin{lemma} \label{stratified ranks characterize forking}
For $a\in G$, $a\thind_A b$ if and only if for each $\Phi, \Theta, k$, we have  that $\thorn_{\Phi, \Theta, k}^G(\tp(a/Ab))= \thorn_{\Phi, \Theta, k}^G(\tp(a/A))$.  \end{lemma}
\noindent
{\em Proof.}
For the left to right direction, if $a\thind_A b$, then all the local ranks of $\tp(a/Ab)$ and $\tp(a/A)$ are equal, and in particular all of the stratified ranks are equal.  

For the right to left direction, first suppose that $\tp(a/Ab)$ \th-divides over $A$. This is witnessed by some $\varphi(x,b)$ and $\theta(y,c)$. Thus one has $\{\varphi^*(x,e\tilde b):e\tilde b\models\widetilde\theta(x,ec)\}$ is $k$-inconsistent.  Hence, 
$$\thorn_{\varphi, \theta, k}^G(\tp(a/A))>\thorn_{\varphi, \theta, k}^G(\tp(a/A)\cup\{\varphi^*(x,eb)\})\geq\thorn_{\varphi, \theta, k}^G(\tp(a/Ab)).$$

Now suppose that $a\nthind_A b$.  Thus for some $n$, $\tp(a/Ab)$ implies a disjunction of formulas $\varphi_i(x,e_i)$, $i\leq n$, such that each $\thorn$-divides over $A$ (witnessed, say, by $\theta_i$ and $k_i$ respectively).  As only the type of the $e_i$ over $Ab$ matters, we may assume that $a\thind_{Ab}e_1$.  Thus by the left to right direction of this lemma, 
$$\thorn^G_{\varphi_1, \theta_1, k_1}(\tp(a/Abe_1))=\thorn^G_{\varphi_1, \theta_1, k_1}(\tp(a/Ab)).$$  But $\varphi_1, \theta_1$, and $k_1$ witness that $\tp(a/Abe_1)$ \th-divides over $A$, so by the previous paragraph,
$$\thorn^G_{\varphi_1, \theta_1, k}(\tp(a/A))>\thorn^G_{\varphi_1, \theta_1, k_1}(\tp(a/Abe_1)).$$
Thus \th-forking implies that not all of the stratified \th-ranks can be equal. \hfill $\blacksquare$\\

Now it is a simple matter to show the existence of \th-generic types.

\begin{theorem} 
There is a $\thorn$-generic type for $G$ over $A$.
\end{theorem}
\noindent
{\em Proof.}
    This is the same proof as the existence of generic types in simple theories in Proposition 4.1.7 in \cite{simple theories}, after one replaces  $D^*(\pi_i, \varphi_i, k_i)$ with $\thorn^G_{\Phi_i,\Theta_i, k_i}(\pi_i)$. \hfill $\blacksquare$\\

We say that a partial type (or a set defined by this type) is $*$-generic if it can be extended to a complete $*$-generic type.

\begin{proposition}
    Let $\pi(x,A)$ be a partial type extending  $G(x)$.  The the following are equivalent:
    \begin{enumerate}
        \item $\pi$ is $\thorn$-generic for $G$,
        \item $\thorn^G_{\Phi, \Theta, k}(\pi)$ is the maximal possible among types in $S_G(A)$, for all $\Phi, \Theta$ and all $k$,
        \item For any $g\in G$ the partial type $g\pi$ does not \th-fork over $\emptyset$, 
        \item For any $g\in G$ the partial type $g\pi$ does not \th-fork over $A$.
    \end{enumerate}

\end{proposition}
\noindent
{\em Proof.}
    The proofs $(1)\Rightarrow (2)$, $(3)\Rightarrow (4)$, and $(4)\Rightarrow (1)$ are the same as those given in Lemma 4.1.9 of \cite{simple theories}.  We consider $(2)\Rightarrow (3)$.

     Suppose that $g\pi$ \th-forks over $\emptyset$.  Then for some $n$, there are $\varphi_i(x,b_i)$, $i\leq n$, whose disjunction is implied by $g\pi$, and such that each $\thorn$-divides over $\emptyset$. Thus there are $\theta_i, k_i$ such that $\{\varphi_i(x,\tilde b_i):\tilde b_i\models\theta_i(y,c_i)\}$ is $k_i$-inconsistent for each $i$.  Let $\Phi=\{\varphi_1,\dots, \varphi_n\}$ and let $\Theta=\{\theta_1,\dots, \theta_n\}$ and let $k$ be the max of $k_1,\dots, k_n$.  Thus for each extension of $g\pi$ to $b_1,\dots, b_n$ one has (as in the proof of Lemma \ref{stratified ranks characterize forking}) that the $\thorn^G_{\Phi,\Theta,k}$-rank of the extension is less than that of  $G(x)$.  But since for a given $\Phi, \Theta, k$ there is always some extension of the same $\thorn^G_{\Phi,\Theta,k}$-rank to any set, it must be that the $\thorn^G_{\Phi,\Theta,k}$-rank of $g\pi$ is less than that of  $G(x)$.  As the stratified ranks are translation invariant, this contradicts $(2)$. \hfill $\blacksquare$

\begin{remark}
The above was rather easier than the corresponding proofs in the case of simple theories because we restrict our attention to definable groups. In a similar fashion, one could define stratified ranks also for type-definable groups, and prove all the above results in this wider context.  
\end{remark}

\begin{question}\label{existence of *-generics}
We have obtained \th-generics by analyzing the local ranks.  If one instead works with an abstract independence relation $\starind$ without associated local ranks, do $*$-generic types still exist?
\end{question}

As mentioned earlier, one may always find a type whose $\uth$-rank is equal to that of $G$.  In fact, this is true for any independence relation $\starind$ for which $*$-generics exist.

\begin{remark}\label{generics have maximal U-rank}
If $p$ is a $*$-generic type of $G$, then $U^*(p)=U^*(G)$. If $U^*(G)<\infty$ and $p$ is such that $U^*(p)=U^*(G)$, then $p$ is $*$-generic. In particular, since a $\thorn$-generic type exists, $\uth(G)=\uth(p)$ where $p$ is a complete $\thorn$-generic type.
\end{remark}  
\noindent
{\em Proof.}
The proof that every $*$-generic type has maximal $U^*$-rank is the same as the proof of a similar fact for $SU$-rank (e.g. in the remarks at the beginning of Section 5.4 of \cite{simple theories}).  For the other direction, consider any type $p \in S(A)$ extending $G(x)$ of maximal $U^*$-rank. Take $a \models p$ and $b\in G$ such that $a \starind_A b$. Then $U^*(p)=U^*(a/A)=U^*(a/Ab)=U^*(ba/Ab) \leq U^*(ba) \leq U^*(p)$. So $ba \starind Ab$, i.e. $p$ is $*$-generic. \hfill $\blacksquare$\\

It is easy to see that our two notions of generics are not the same: 

\begin{example}
    Let $R$ be an $\omega$-saturated real closed field. We may think of $R^2$ as $C$, the algebraic closure of $R$.  Consider the unit circle in $R^2$ as a multiplicative subgroup of $C$.  Then consider the circle intersect an infinitesimal neighborhood of the point $(1,0)$.  This is of $\uth$-rank one as is $S^1$.  Thus it is \th-generic, but clearly it is not generic.
\end{example}

\noindent On the other hand, we do have the converse:

\begin{proposition}\label{generic is thorn generic}
    A generic type is \th-generic.
\end{proposition}
\noindent
{\em Proof.}
    Let $X\subseteq G$ be a generic definable set.  We must show that each stratified rank of  $X$ is maximal.  Say $G\subseteq g_1 X\cup\dots\cup g_n X$.  Thus  $G(x)$ implies $\bigvee_{i\leq n} (x\in g_i X)$.  Thus the $\thorn^G_{\Phi,\Theta, k}$-rank of $G$ is equal to the maximum of the $\thorn^G_{\Phi,\Theta, k}$-ranks of $g_i X$.  But these are all the same, as the stratified ranks are translation invariant. \hfill $\blacksquare$\\

As in Question \ref{existence of *-generics}, it is not clear whether every generic type is $*$-generic.
One can say slightly more about the relation between generic and \th-generic in the case of subgroups.  The following is clear, but useful.

\begin{proposition} \label{bounded index is generic}
1. Let $H$ be a type definable subgroup of $G$.   $H$ has bounded index in $G$ $\Leftrightarrow$ $H$ is generic in $G$ $\Rightarrow$ $H$ is \th-generic in $G$.\\
2. Let $H$ be a definable subgroup of $G$.  $H$ has finite index in $G$ $\Leftrightarrow$ $H$ is generic in $G$ $\Leftrightarrow$ $H$ is \th-generic in $G$.
\end{proposition}

In the following considerations, we will often use cosets modulo a definable subgroup $H$ of $G$. Then, for every $g \in G$, $\bar g$ will always denote the coset $gH$ treated as an element of $\C^{eq}$.

The following will be useful later.
\begin{proposition}\label{quotients}
Let $H$ be an $A$-definable normal subgroup of $G$.  If $\tp(g/A)$ is a $*$-generic of $G$, then $\tp(\bar g/A)$ is a $*$-generic of $G/H$. 
Furthermore, if at least one $*$-generic type of $G$ exists (e.g. it is the case when $*=\thorn$), then all $*$-generics of $G/H$ over $A$ arise in this fashion.
\end{proposition}
\noindent
{\em Proof.} Assume for simplicity that $A=\emptyset$. Take any $\bar h\in G/H$ such that $\bar g \starind \bar h$. We need to show that $\bar h \bar g \starind \bar h$. 
First we can find $g' \models \tp(g/\bar g)$ so that $g' \starind_{\bar g} \bar h$. Then $g' \starind \bar h$. Now we choose $h' \models \tp(h/\bar h)$ so that $h' \starind_{ \bar h} g'$. So $g' \starind h'$. Since $\tp(g')=\tp(g)$ is $*$-generic, we get $h'g' \starind h'$. As $\bar h=\bar h' \in \acl(h')$ and $\bar g=\bar g'$ imply $\bar h \bar g = \bar h' \bar g' =\barxx{h'g'} \in \acl(h'g')$, we get $\bar h \bar g \starind \bar h$.

For the converse, take any $\bar g$ which is $*$-generic in $G/H$. Now choose $h \in G$ which is $*$-generic in $G$ and $h \starind g$. Then $\bar g \starind h$, so 
$\bar g$ is $*$-generic of $G/H$ over $h$. Since $\bar{h} \in \acl(h)$, we get  
$\bar h \bar g \starind h$. Choose $g' \models \tp(hg/\bar h \bar g)$ so that $g' \starind_{\bar h \bar g} h$. Then $h \starind g'$. Since $\tp(h)$ is $*$-generic, we get that $g_1:=h^{-1}g '$ is also $*$-generic. On the other hand, $\bar {g'}=\barxv{hg}$ so $\bar g_1=\barxxx{h^{-1}g'}=\barxxxv{h^{-1}hg}=\bar g$. \hfill $\blacksquare$

\begin{definition}
Suppose $H$ is an $A$-definable subgroup of $G$ and a coset $aH$ is definable over $B$. A type $p \in S(B)$ is a \emph{$*$-generic of $aH$ over $B$} if there is a non-$*$-forking extension $q \in S(ABa)$ of $p$ such that $a^{-1}q$ is $*$-generic of $H$.
\end{definition}

It is obvious that $p$ is $\thorn$-generic for $aH$ if and only if it is a type extending $x \in aH$ of maximal possible stratified ranks (equal to the stratified ranks of $H$). If   $U^*(G)<\infty$, then $p$ is $*$-generic of $aH$ if and only if it is a type extending $x \in aH$ of maximal possible $U^*$-rank (and $U^*(p)=U^*(H)$). One can also apply the proof of \cite[Lemma 4.3.12]{simple theories} to conclude that if $g$ is $\thorn$-generic for $G$ over $A$, then $gH$ is definable over $A,\bar g$ and $g$ is $\thorn$-generic for $gH$ over $A, \bar g$ (the proof uses the existence of a $\thorn$-generic in $H$). 

In the final part of this section, we check some basic properties of $U^*$-rank in groups. In particular, we show a version of Lascar inequalities for groups.


\begin{proposition}\label{quotients1}
Assume $H$ is an $A$-definable subgroup of $G$. If $\tp(g/A)$ is a $*$-generic type of $G$, then $U^*(G/H)=U^*(\bar g/A)$.
Furthermore, if  $U^*(G)<\infty$ and at least one $*$-generic type of $G$ exists, then all elements of maximal $U^*$-rank in $G/H$ over $A$ arise in this way.
\end{proposition}
\noindent
{\em Proof.} Assume for simplicity that  $A= \emptyset$. Take any $h \in G$. We need to show that $U^*(\bar g) \geq U^*(\bar h)$. Wlog $h \starind g$. Then $\bar g \starind h$. Since $g$ is $*$-generic, we get that $hg\starind h$. Let $g_1= hgh^{-1}$. Since $g_1h=hg$, we see that $g_1h \starind h$ and hence $\barxx{g_1h} \starind h$. Using this, we get $U^*(\bar g)=U^*(\bar g/h)=U^*(\barxxxv{h^{-1}g_1h}/h)=U^*(\barxv{g_1h}/h)=U^{*}(\barxv{g_1h}) \geq U^*(\barxx{g_1h}/g_1)=U^*(\bar h/g_1)\geq U^*(\bar h)$.

The proof of the second part is similar as in Proposition \ref{quotients}. Suppose that $\bar g\in G/H$ is such that $U^*(\bar g)=U^*(G/H)$. Now choose $h \in G$ which is $*$-generic in $G$ and $h \starind g$. Then $\bar g \starind h$, so $U^*(\bar g/h)=U^*(G/H)$. Since $\bar g$ and $\barx{hg}$ are interalgebraic over $h$, we get that $U^*(\barx{hg}/h)=U^*(G/H)\geq U^*(\barx{hg})$. So $\barx{hg} \starind h$. Choose $g' \models \tp(hg/\barx{hg})$  so that $g' \starind_{\barx{hg}} h$. Then $h \starind g'$. Since  $\tp(h)$ is $*$-generic, we get that $g_1:=h^{-1}g '$ is also $*$-generic. On the other hand, $\bar{g'}=\barx{hg}$ so $\bar g_1=\barxxx{h^{-1}g'}=\barxxxv{h^{-1}hg}=\bar g$. \hfill $\blacksquare$

\begin{proposition}[Lascar inequalities for groups]
Let $H$ be a definable subgroup of $G$. Then\\
1. $\alpha+U^*(G/H) \leq U^*(G) \leq U^*(H) \oplus U^*(G/H)$ for every $\alpha<U^*(H)$ or $\alpha=U^*(H)$ if there is a $*$-generic type of $H$,\\
2. $\uth(H)+\uth(G/H)\leq\uth(G)\leq \uth(H)\oplus\uth(G/H)$.
\end{proposition}
\noindent
{\em Proof.} (2) follows from (1) and the existence of $\thorn$-generics.\\
(1). Wlog $H$ is $\emptyset$-definable. First let us prove the right inequality. Take any $\lambda < U^*(G)$ or $\lambda=U^*(G)$ if there is a $*$-generic type of $G$. Then, by Remark \ref{generics have maximal U-rank}, there is $g \in G$ such that $U^*(g) \geq \lambda$. In the following computation, $gH$ denotes a subset of $G$ and $\bar g=gH$ an element of $\C^{eq}$. Using Lascar inequalities, we get $\lambda \leq U^*(g) =U^*(\bar g/g)+U^*(g) \leq U^*(g,\bar g) \leq   U^*(g/\bar g) \oplus U^*(\bar g) \leq U^*(gH)\oplus U^*(G/H)=U^*(H)\oplus U^*(G/H)$.

Now we turn to the left inequality. Take any $\alpha$ as in the proposition and $\beta<U^*(G/H)$ or $\beta=U^*(G/H)$  there is an element in $G/H$ of maximal $U^*$-rank. It is enough to show that $\alpha+\beta \leq U^*(G)$. By Remark \ref{generics have maximal U-rank}, there are $h \in H$ and $\bar g\in G/H$ such that $U^*(h) \geq \alpha$ and $U^*(\bar g) \geq \beta$. Wlog we can assume that $h \starind g$. Then $U^*(gh/\bar g) \geq U^*(gh/g)=U^*(h/g)=U^*(h) \geq \alpha$. Let $g_1=gh$. Then $\bar g=  \bar g_1$ and so, by Lascar inequalities, we get $U^*(G) \geq U^*(g_1)=U^*(\bar g_1/g_1)\oplus U^*(g_1) \geq U^*(g_1,\bar g_1) \geq U^*(g_1/\bar g_1) + U^*(\bar g_1)=U^*(gh/\bar g) + U^*(\bar g) \geq \alpha + \beta$. \hfill $\blacksquare$

\subsection*{\th-orthogonality and \th-regular types}
We work in our rosy theory $T$. We define all notions with respect to an arbitrary independence relation $\starind$. In particular, everything applies to $\thorn$-independence.
We define $*$-orthogonality and $*$-regular types in the same way as the corresponding notions are defined in stable theories.

\begin{definition} 
Let $p$ and $q$ be complete types, and $A$ be a set containing each of their domains.  If $a\starind_A b$ for any $a,b$ realizing non-$*$-forking extensions to $A$ of $p,q$ respectively, then we say that $p$ and $q$ are \emph{$*$-orthogonal}.  We say that $p \in S(A)$ is \emph{$*$-regular} if it is $*$-orthogonal to all its $*$-forking extensions.
\end{definition}

\begin{remark}
$*$-regularity is preserved under non-$*$-forking extensions.
\end{remark}

As in the stable case, using Lascar inequalities one can show the following:

\begin{remark}\label{omega to the alpha implies regular}
Each type of $U^*$-rank $\omega^{\alpha}$ is $*$-regular.
\end{remark}

Now we study $*$-regularity in our group $G$ definable in $\C$.

\begin{lemma}\label{generic types translated}
    Let $p,q \in S(A)$ extend $G(x)$ and let $p$ be $*$-generic.  Then there is some $g$ such that $gp\cup q$ is a non-$*$-forking extension of $q$  (i.e. there is a complete type extending $gp\cup q$ which does not $*$-fork over $A$). 
\end{lemma}

\noindent \emph{Proof.}
     Take $a\models p$ and $b\starind_A a$ with $b\models q$. Since $p$ is $*$-generic, so is $\tp(a^{-1}/A)$. Since $a$ and $a^{-1}$ are interdefinable, we also have $a^{-1}\starind_A b$. Thus we see that $ba^{-1}\starind Ab$. We note that $b\models ba^{-1}p\cup q$, a partial type over $A(ba^{-1})$.  Since $b\starind_A ba^{-1}$, we may extend this partial type to a complete type, $q'$, over $A(ba^{-1})$ extending $q$ which does not $*$-fork over $A$.  We let $g:=ba^{-1}$. \hfill $\blacksquare$

%

\begin{proposition}\label{regular generic orthogonal to non thorn generic}
If $p \in S(A)$ is a $*$-regular $*$-generic type of $G$, then it is $*$-orthogonal to every non-$*$-generic type $q \in S(A)$ of $G$.
\end{proposition}
\noindent
{\em Proof.} 
%
Wlog $A= \emptyset$. 
By Lemma \ref{generic types translated}, we can choose $g$ and a non-$*$-forking extension $r \in S(B)$ of $q$ so that $gp \subseteq r$, where $B=\dcl(g)$. Then, by Lemma \ref{basic facts about *-generics}.2 and \ref{basic facts about *-generics}.3, we get $g^{-1}r \in S(B)$ is a $*$-forking extension of $p$, so by $*$-regularity of $p$, $g^{-1}r$ is $*$-orthogonal to $p$. Hence 

\begin{enumerate}
\item[($*$)] $r$ is $*$-orthogonal to $p$.
\end{enumerate} 

Now suppose for a contradiction that $p$ and $q$ are not $*$-orthogonal. Then  there is $C\subseteq G$, $a \models p$, and $b \models q$ such that $a \starind C$, $b \starind C$, and $a \nstarind_C b$.

Take any $b'\models r$. Since $\tp(b)=\tp(b')$, we can choose $g' \models \tp(g)$ so that $\tp(gb')=\tp(g'b)$ and $g' \starind_{b} Ca$. As $g \starind b'$, we get $g' \starind b$ and hence $g' \starind Cab$. Thus $a \starind Cg'$, $b \starind Cg'$, and $a \nstarind_{Cg'}b$.  It follows that $\tp(b/g')$ is not $*$-orthogonal to $p$. This yields a contradiction with ($*$). \hfill $\blacksquare$

\section{Groups of $\uth$-rank 1}

In this section $G$ is a definable group in a monster model $\C$ of a rosy theory $T$. In all results of this section in which we assume that $T$ satisfies NIP, one can replace this assumption by the weaker condition that $G$ has icc on centralizers (see Corollary \ref{icc on centralizers}).

\begin{proposition}\label{thorn generic involution}
If $G$ contains a $\thorn$-generic involution, then it contains a $\thorn$-generic element $g$ such that $[G:C(g)]<\omega$.
\end{proposition}
\noindent {\em Proof.} Let $i$ be a $\thorn$-generic involution. Choose $j \models \tp(i)$ so that $j \thind i$. Then $j$ is $\thorn$-generic over $i$. So $i \thind ij$. Thus $\tp(i/ij)$ is $\thorn$-generic.

On the other hand, $(ij)^i=(ij)^{-1}$. So for every $i' \models \tp(i/ij)$ we
also have $(ij)^{i'} =(ij)^{-1}$. Hence $(ij)^i=(ij)^{i'}$, so $i'\in C(ij)i$. We
conclude that $C(ij)i$ is definable over $ij$ and a formula defining this set
belongs to $\tp(i/ij)$. Since $\tp(i/ij)$ is $\thorn$-generic, we get that $C(ij)i$ and hence $C(ij)$ is $\thorn$-generic.
Hence $[G:C(ij)] < \omega$. Moreover, as $i$ is $\thorn$-generic and $i \thind j$, we get that $ij$ is also $\thorn$-generic. \hfill $\blacksquare$

\begin{corollary}\label{involution corollary 1}
If $T$ satisfies NIP and $G$ has a $\thorn$-generic involution, then $G$ is abelian-by-finite.
\end{corollary} 

\noindent {\em Proof.}

Let $H$ be the centralizer connected component of $G$ (a definable set by the icc on centralizers).  It is of finite index, and we will show that it is abelian.  Let $g$ be a \th-generic involution; by Proposition \ref{thorn generic involution}, $C(g)$ is a finite index subgroup of $G$.  Thus $H$ is a subgroup of $C(g)$, and hence $g\in C_G(H)$.  Since $C_G(H)$ may be extended to a \th-generic type, it is of finite index in $G$, and thus $C_G(H)=H$, and the latter is abelian. \hfill $\blacksquare$\\

The next two corollaries follow immediately from Proposition \ref{thorn generic involution} and Corollary \ref{involution corollary 1}. 

\begin{corollary}\label{involution corollary 2}
Assume $\uth(G)=\alpha<\infty$. If $G$ contains an involution of $\uth$-rank $\alpha$, then it contains an element $g$ such that $[G:C(g)]<\omega$ and $\uth(g)=\alpha$.
\end{corollary}

\begin{corollary}\label{involution corollary 3}
If $T$ satisfies NIP, and $G$ has an involution $i$ such that $\uth(i)=\uth(G)<\infty$, then $G$ is abelian-by-finite.
\end{corollary} 

\begin{theorem}\label{definable abelian subgroup}
If $T$ satisfies NIP, $G$  has hereditarily fsg, and $0<\uth(G)<\infty$, then $G$ contains an infinite definable abelian subgroup.
\end{theorem}
\noindent {\em Proof.} 
We can replace $G$ by an infinite definable subgroup of least possible $\uth$-rank and we can assume that $G$ is centralizer connected. We will show that $G$ is abelian.

If $Z(G)$ is infinite, then $[G:Z(G)] < \omega$, so $G=Z(G)$, and we are done. So we can assume that $Z(G)$ is finite. Let $H=G/Z(G)$. Then $H$ is infinite and each non-trivial element of $H$ has a finite centralizer.  Now we will show that this leads to a contradiction.\\[2mm]
{\bf Claim}  {\em There are finitely many conjugacy classes in $H$.}\\[2mm]
\noindent {\em Proof.} Take any $g \in H\setminus \{ e \}$. Since $C(g)$ is finite, we get that $\uth(g^H)=\uth(H)$. Now the relation of being in the same conjugacy class is a $\emptyset$-definable equivalence relation on $H \setminus \{ e \}$ whose classes are $\thorn$-generic. Hence there must be only finitely many of them. \hfill $\square$\\[2mm]
By the above claim, we get that $H^{00}$ is definable. So wlog $H=H^{00}$. By Proposition \ref{unique generic}, we get that there is a unique generic type in $H$. So in virtue of the claim, we get that there is a unique generic conjugacy class $a^H$.\\[2mm]
{\bf Case 1.} {\em There is no involution in $H$.}

Since $a^H$ is generic, $(a^{-1})^{H}$ is also generic. By uniqueness, $a^H=(a^{-1})^H$. Hence there is $g \in H$ such that $a^{-1}=g^{-1}ag$. Thus $a=g^{-2}ag^2$. Since $a$ is not an involution, we get $a \in C(g^2) \setminus C(g)$. So $C(g) \subsetneq C(g^2)$. Since $C(g)$ is finite and there are no involutions, we get that $g$ has an odd exponent. This implies $C(g^2)= C(g)$, a contradiction.\\[2mm]
{\bf Case 2.} {\em There is an involution $i \in H$.}

There are two ways to get a contradiction. By \cite[Theorem 2.1]{locally finite groups}, we get that there is a non-trivial element with infinite centralizer, a contradiction. Alternatively, we can argue as follows: there is $g \in H$ such that $i^g$ is $\thorn$-generic, but $i^g$ is also an involution, so by Proposition \ref{thorn generic involution} (or Corollary \ref{involution corollary 1}), we get a non-trivial element with an infinite centralizer, a contradiction. \hfill $\blacksquare$\\

If $G$ has fsg, then every definable subgroup of finite index in $G$ also has fsg. 
So if every definable subgroup of $G$ is either finite or of finite index in $G$ and $G$ has fsg, then $G$ has hereditarily fsg.
Thus we have the following immediate corollary:

\begin{corollary}\label{Cor 2.6}
If $T$ satisfies NIP, $G$ has fsg, $\uth(G)<\infty$, and each definable subgroup of $G$ is either finite or of finite index in $G$, then $G$ is abelian-by-finite.
\end{corollary} 

In particular, we have proven Theorem \ref{Main Theorem 1}:

\begin{proven main theorem}
If T has NIP, $G$ has fsg, and $\uth(G)=1$, then $G$ is abelian-by-finite.
\end{proven main theorem}

Now we are going to modify the proof of Theorem \ref{definable abelian subgroup} (using the idea of the proof of \cite[Proposition 7.2]{groupes stables}) to generalize Theorem \ref{Main Theorem 1} to the case of groups with a \th-regular \th-generic type.

\begin{theorem} \label{regular generic implies abelian by finite}
If $T$ satisfies NIP, $G$ has fsg, and at least one \th-regular \th-generic type, then $G$ is abelian-by-finite.
\end{theorem}
\noindent {\em Proof.} First of all we can assume that $G$ is centralizer connected.\\[2mm]
{\bf Claim 1} {\em The conjugacy class of every non-central element is $\thorn$-generic ( i.e. there is an element in this conjugacy class which is $\thorn$-generic over a name of this class).} \\[2mm]
\noindent {\em Proof.} Take any $a \in G$ such that $a^G$ is not $\thorn$-generic. Let  $\tp(b/a)$ be a \th-regular \th-generic type. Then $b^{-1}ab$ is not $\thorn$-generic over $a$. By Proposition \ref{regular generic orthogonal to non thorn generic}, we get $b \thind_a b^{-1}ab$. Thus $b$ is $\thorn$-generic over $\{a, b^{-1}ab\}$. So the set defined by the formula $x^{-1}ax=b^{-1}ab$ is $\thorn$-generic and definable over $\{a, b^{-1}ab\}$. But this set is equal to $C(a)b$. Hence $C(a)$ is $\thorn$-generic. Thus $[G:C(a)] < \omega$. Since $G$ is centralizer connected, we get $a \in Z(G)$. \hfill $\square$\\[2mm]
{\bf Claim 2} {\em $G$ is the union of $Z(G)$ and finitely many $\thorn$-generic conjugacy classes.}\\[2mm]
\noindent {\em Proof.} By Claim 1, $G$ is the union of $Z(G)$ and some number of $\thorn$-generic conjugacy classes. Suppose for a contradiction that there are infinitely many of them. Then the relation of being in the same conjugacy class is $\emptyset$-definable and it divides $G\setminus Z(G)$ into infinitely many classes. So at least one of these classes, say $C$, is non-algebraic over $\emptyset$ (as an element of $\C^{eq}$). Now we can choose $c \in C$ which is $\thorn$-generic over $C$. Since $C \in \dcl(c)$ and $C \notin \acl(\emptyset)$, we get $c \nthind C$, a contradiction with Lemma \ref{basic facts about *-generics}.1. \hspace*{416pt} $\square$\\[2mm]
We will show that $Z(G)=G$. Suppose for a contradiction that it is false. Then $[G:Z(G)] \geq \omega$. By Claim 2 and Proposition \ref{quotients}, $G/Z(G)$ is the union of $\{ e \}$ and finitely many $\thorn$-generic conjugacy classes.  Hence $H:=(G/Z(G))^{00}$ is definable and it is also the union of $\{ e \}$ and finitely many $\thorn$-generic conjugacy classes. 
In the same way as in the proof of Theorem \ref{definable abelian subgroup}, we conclude that there is a unique generic conjugacy class $a^H$. 

As in the proof of Theorem \ref{definable abelian subgroup}, we get that there is an involution $i \in H$. Since $i^H$ is $\thorn$-generic, we get that there is an involution which is $\thorn$-generic. We finish using Proposition \ref{thorn generic involution}. \hfill $\blacksquare$\\  

By Theorem \ref{regular generic implies abelian by finite} and Remark \ref{omega to the alpha implies regular}, we get the following strengthening of Theorem \ref{Main Theorem 1}, which generalizes the appropriate result about superstable groups.

\begin{corollary}\label{Cor 2.14}
If $T$ satisfies NIP, $G$ has fsg, and $\uth(G)=\omega^{\alpha}$, then $G$ is abelian-by-finite.
\end{corollary}

At the end let us make a few remarks. Proposition \ref{thorn generic involution} and Corollaries \ref{involution corollary 1}, \ref{involution corollary 2},  and \ref{involution corollary 3} are true for an arbitrary independence relation,  $\starind$ (the same proofs work). Since $\uth$-rank is less or equal than the $U^{*}$-rank, Theorem \ref{definable abelian subgroup}, Corollary \ref{Cor 2.6} and Theorem \ref{Main Theorem 1} are also true for $\ind^*$. Since Proposition \ref{regular generic orthogonal to non thorn generic} is true for $\starind$, one can easily check that the whole proof of Theorem \ref{regular generic implies abelian by finite} also works for $\starind$.

As to Corollary \ref{Cor 2.14}, it is true for an arbitrary independence relation $\ind^*$ under the additional assumption that there is a type of $U^{*}$-rank $\omega^{\alpha}$ (in other words a $*$-generic exists), in which case it follows immediately from the fact that Theorem \ref{regular generic implies abelian by finite} is true for $\starind$; an alternative way to prove this is to modify the proof of Theorem \ref{definable abelian subgroup} and use Lascar inequalities for groups.

By Proposition \ref{generic is thorn generic} and Theorem \ref{regular generic implies abelian by finite}, we get that if $T$ satisfies NIP, $G$ has fsg and at least one \th-regular generic type, then $G$ is abelian-by-finite.

Notice that the proofs of Theorems \ref{Main Theorem 1} and \ref{regular generic implies abelian by finite} produce a definable abelian subgroup of finite index in $G$, namely the centralizer connected component of $G$. 

Notice also that when $A$ is an arbitrary abelian subgroup of finite index in $G$, we may follow the proof of Theorem 3.17 of \cite{groupes stables} in considering the intersection of $C(a)$ for $a\in A$.  This is definable, as is its center, which is abelian of finite index in $G$ and contains $A$.
 
Finally, we should mention a conjecture generalizing the results above.  The existence of finitely satisfiable generics implies the weaker condition that $G$ is \emph{definably amenable}. (This means, roughly, that there is a left invariant probability measure on the definable sets of $G$.  See \cite{NIP paper} for a precise definition.)

\begin{conjecture} 
In each result in this section, the hypothesis ``$G$ has fsg'' may be replaced with ``$G$ is definably amenable''.
\end{conjecture}

\section{Groups of $\uth$-rank 2}

In this section $G$ is a group definable in a monster model $\C$ of an arbitrary theory $T$. 

\begin{theorem} \label{3.1}
If $G$ has hereditarily fsg, the definable quotients of definable subgroups of $G$ satisfy icc on centralizers, and $\uth(G)$=2, then $G$ is solvable-by-finite.
\end{theorem}

By Corollary \ref{icc on centralizers}, Theorem \ref{3.1} implies Theorem \ref{Main Theorem 2}: 

\begin{proven main theorem}
If $T$ satisfies NIP, $G$ has hereditarily fsg, and $\uth(G)=2$, then $G$ is solvable-by-finite.
\end{proven main theorem}
The rest of this section is devoted to the proof of Theorem \ref{3.1}. 
Some ideas are taken from the proof that each connected group of Morley rank 2 is solvable \cite[Theorem 3.16]{groupes stables}. The main obstacle in comparison with the Morley rank 2 case is that a $\thorn$-generic may not be generic and there may be many $\thorn$-generics even in the connected component. 

The structure of the proof is as follows. We suppose for a contradiction that $G$ is not solvable-by-finite. First we define Borels (albeit in a slightly different manner from that which is used in the Morley rank 2 case because we do not have definable connected components), and we study their properties. Then we use them to find involutions. In the last part of  the proof we use Borels, involutions and some particular function that comes from the theory of black box groups to get a final contradiction.\\[2mm] 
{\em Proof of Theorem  \ref{3.1}.}
By icc, we can assume that $G$ is centralizer connected. If $Z(G)$ is infinite, then either $[G:Z(G)] <\omega$ or $\uth(Z(G))=\uth(G/Z(G))=1$ and then we are done by Theorem \ref{Main Theorem 1}. So we can assume that $Z(G)$ is finite.  Then $G/Z(G)$ is centerless  and centralizer connected. So wlog $G$ is centerless. Suppose for a contradiction that $G$ is not solvable-by-finite.

If centralizers of all non-trivial elements in $G$ are finite, we can argue in the same way as in the proof of Theorem \ref{definable abelian subgroup}. So we can assume that there is a non-trivial element in $G$ with an infinite centralizer.

\begin{definition}
We say that a subgroup $B$ of $G$ is a \emph{Borel} if it is a minimal infinite intersection of centralizers.
\end{definition}

By icc and the last paragraph, we have that at least one Borel exists and all Borels are intersections of finitely many infinite centralizers, so they are definable.
\begin{claim} (i) Every Borel has infinite index in $G$.\\
(ii) All Borels are abelian.\\
(iii) Any two Borels are either equal or they have trivial intersection.\\
(iv) Every Borel has finite index in its normalizer.\\
(v) Every non-trivial element $a \in G$ with an infinite centralizer centralizes exactly one Borel, which will be denoted by $B(a)$. In fact, $B(a)$ is the centralizer connected component of $C(a)$ computed in $G$.
\end{claim}
\noindent {\em Proof.} (i) This follows from the fact that the centralizer of any non-trivial element has infinite index in $G$.\\
(ii) By (i), $\uth(B)=1$. So by Theorem \ref{Main Theorem 1}, there is an abelian subgroup $B_0$ of finite index in $B$. By the definition of Borels, $B$ is centralizer connected. Hence $B_0 \leq Z(B)$. So for every $a \in B$, $[B:C_B(a)]<\omega$. Since $B$ is centralizer connected, it must be abelian.\\
(iii) Suppose $B_1 \ne B_2$ are Borels. Then $B_1 \cap B_2$ is finite. Let $a \in B_1 \cap B_2$. Then $B_1,B_2 \subseteq C(a)$. So $[C(a): B_1]\geq \omega$, which implies $\uth(C(a))=2$, so $[G:C(a)] < \omega$. Hence $a=e$.\\
(iv) If $B$ is a Borel and $[N(B):B]\geq \omega$, then $\uth(N(B)/B)=1$ and $[G:N(B)] < \omega$, so we finish using Theorem \ref{Main Theorem 1}.\\
(v) Since $\uth(C(a))=1$, if $a$ centralizes a Borel $B$, then $B$ is a subgroup of finite index in $C(a)$. So by (iii), there is a unique Borel centralized by $a$, and, of course, it must be the centralizer connected component of $C(a)$ computed in $G$. \hfill $\square$\\

The following is an easy corollary of Claim 1.

\begin{claim} If $B_1 \ne B_2$ are Borels, then:\\ 
(i) $C(B_1) \cap C(B_2) = \{ e \}$,\\
(ii) $N(B_1) \cap N(B_2)$ is finite.
\end{claim}
\noindent {\em Proof.} The first part follows from the fact that every non-trivial element centralizes at most one Borel. The second part is an easy consequence of  the fact that   $B_1 \cap B_2 = \{ e \}$ and $B_1$ and  $B_2$ are subgroups of finite index in $N(B_1)$ and $N(B_2)$, respectively. \hfill $\square$ 

\begin{claim}
If $B$ is a Borel and $a$ is generic over a name of $B$, then $BaB$ is generic.
\end{claim}
\noindent {\em Proof.} 
We have that there are $m,n \in \omega$ such that:
\begin{enumerate} 
\item[(1)] $BN(B)B=Bb_1B \cup \dots \cup Bb_mB$ for some $b_1,\dots, b_m \in N(B)$, 
\item[(2)] $B(G \setminus N(B))B= Ba_1B \cup \dots \cup Ba_nB$ for some $a_1, \dots, a_n \in G \setminus N(B)$.
\end{enumerate}

Item (1) follows from Claim 1(iv).
To see $(2)$, notice that if $a \notin N(B)$, then $B^a \ne B$, hence $B^a \cap B
= \{e\}$. So $f:B \times B \to BaB$ defined by $f(x,y)=xay$ is a definable
bijection. Thus $\uth(BaB)=2$. Moreover, for all $a_1,a_2 \in G$ we have that $Ba_1B$ and $Ba_2B$ are either equal or disjoint. Since $\uth(G)=2$, we
get that there are only finitely many sets of the form $BaB$ for $a \notin N(B)$.

Let $d$ be a name for $B$. Take any generic $a$ over $d$. Let $p=\tp(a/d)$. 
By (1) and (2), we get that $p(G)$ is covered by $Ba_1'B, \dots, Ba_k'B$ for some $a_1',\dots, a_k'\in p(G)$

Let ${\bf p'}$ be a global generic extension of $p$. Then there is $i$ such that a formula defining the set $Ba_i'B$ is in ${\bf p'}$. So we have proved that there is $a' \models p$ such that $Ba'B$ is generic. Hence  $BaB$ is generic. \hfill $\square$\\

In the remainder of the proof, we analyze $B\cap G^{00}$.  The next claim is not necessary to finish the proof, but we include it, since it gives us a better understanding of Borels (and lets us simplify our notation).

\begin{claim}
Every Borel is relatively connected in $G^{00}$, i.e. $B^{00}= B \cap G^{00}$.
\end{claim}
\noindent {\em Proof.} The inclusion $B^{00} \subseteq B \cap G^{00}$ is obvious. We will prove $B^{00} \supseteq B \cap G^{00}$. 

Suppose for a contradiction that there is $c \in (B \cap G^{00}) \setminus B^{00}$. Let $a \in G^{00}$ be generic over a name of $B$. 

By Claim 3, we easily get that boundedly many two-sided translates of $B^{00}aB^{00}$ cover $G$, so all these translates are defined by partial generic types. Let $\bf p$ be a global generic type extending a partial generic type defining the set $cB^{00}aB^{00}$.  Then a partial type defining the set $B^{00}aB^{00}$ is contained in $c^{-1}{\bf p}={\bf p}$ (as $c \in G^{00}$ and $Stab({\bf p})=G^{00}$).
Hence we get

$$B^{00}aB^{00} \cap cB^{00}aB^{00} \ne \emptyset.$$

Combining this with the fact that the function $f : B \times B \to G$ defined by $f(b_1,b_2)=b_1ab_2$ is injective and that $cB^{00} \cap B^{00}=\emptyset$, we get a contradiction. \hfill $\square$

\begin{claim} If $B$ is a Borel and $a \in N(B) \setminus C(B)$, then $C(a)$ is finite.
\end{claim}
\noindent {\em Proof.} First notice that we have
\begin{enumerate}
\item[$(*)$] $a^{N(B)}$ is infinite.
\end{enumerate}
Otherwise, $a^B$ would also be finite. So $C(a) \cap B$ would be infinite. Hence $C(a)\cap B =B$, which would mean that $a \in C(B)$, a contradiction.

For any $g$ we have $a^{N(B)g}=(a^{N(B)})^g \subseteq N(B)^g=N(B^g)$. Using this and Claim 2(ii) we see that 
\begin{enumerate}
\item[$(**)$] if $N(B)g_1 \ne N(B)g_2$, then $a^{N(B)g_1} \cap a^{N(B)g_2}$ is finite.
\end{enumerate}
 
Now add $a$ and a name of $B$ to the language. Take $g$ so that $N(B)g \notin \acl^{eq}(\emptyset)$. Let $\phi(x,y)$ be a formula over $\emptyset$ such that $\phi(\C,N(B)g)=a^{N(B)g}$ (on the left hand side $N(B)g$ is treated as an element of the sort $N(B) \backslash G$ and on the right hand side it is treated as a set).  By $(*)$, $(**)$ and Lemma \ref{finite intersections good enough}, we get that there is $b \in \phi(\C,N(B)g)=a^{N(B)g}$ such that $\uth(b)=2$. Hence, working in the original language, $\uth(b/a)=2$. Since $a^G=b^G$, we get that $\uth(a^G)=2$, so $C(a)$ is finite. \hfill $\square$\\ 

Now we will combine the above result and the proof of Claim 4 to get that Borels are relatively self-normalizing in $G^{00}$, i.e for every Borel $B$ we have $N_{G^{00}}(B\cap G^{00})=B \cap G^{00}$.  By Claim 1(iii), this is equivalent to  
the condition $N(B) \cap G^{00} = B \cap G^{00}$, and by Claim 4, to the statement $N_{G^{00}}(B^{00})=B^{00}$ .

\begin{claim} All Borels are relatively self-normalizing in $G^{00}$.
\end{claim}
\noindent {\em Proof.} Take any Borel $B$ and $a \in G^{00}$ generic over a name of $B$. Suppose for a contradiction that there is $c \in (N(B) \cap G^{00}) \setminus (B \cap G^{00})$. The same argument as in the proof of Claim 4 yields
$B^{00}aB^{00} \cap cB^{00}aB^{00} \ne \emptyset$. Hence
$$(!) \;\;\;\;\;\; BaB \cap cBaB \ne \emptyset.$$
{\bf Subclaim} {\em The function $f: N(B) \times B \to G$ defined by $f(b_1,b_2)=b_1ab_2$ is injective.}\\[2mm]
\noindent {\em Proof.} Suppose $b_1ab_2=c_1ac_2$ for some $b_1,c_1 \in N(B)$ and $b_2,c_2 \in B$. Let $b=c_1^{-1}b_1$. We see that $b^a \in B$, so $b \in B^{a^{-1}}$. Hence $C(b)$ is infinite. On the other hand, $b \in N(B)$, so by Claim 5, we get that $b \in C(B)$. 

Since $a \notin N(B)$, we have that $B^{a^{-1}} \ne B$, so by Claim 2(i), $C(B) \cap B^{a^{-1}} = \{ e \}$.  Since $b \in C(B) \cap B^{a^{-1}}$, we get $b=e$. So $b_1=c_1$ and then also $b_2=c_2$. \hfill $\square$\\[2mm]   
Subclaim and $(!)$ yield a contradiction. \hfill $\square$\\ 

So far we have been studying various properties of Borels. Now we will use Borels to find an involution in $G^{00}$.

 \begin{claim}
There is an involution in $G^{00}$.
\end{claim}
\noindent {\em Proof.} By fsg, there is a generic element $a \in G^{00}$. Moreover, we have at least one Borel $B$. So we can choose $a \in G^{00}$ generic over a name of $B$. 

By Claim 3, we get 
\begin{enumerate}
\item[$(*)$] $BaB$ and $Ba^{-1}B$ are generic.
\end{enumerate}
{\bf Subclaim.} $B^{00}aB^{00}=B^{00}a^{-1}B^{00}$.\\[2mm]
\noindent {\em Proof.} It is clear that $B^{00}aB^{00}$ and $B^{00}a^{-1}B^{00}$ are contained in $G^{00}$.  Let $\pi$ be a partial type defining the set $B^{00}a^{-1}B^{00}$. 

Since by $(*)$, boundedly many two-sided translates of $B^{00}a^{-1}B^{00}$ cover $G$, $\pi$ is a partial generic type. Let $\bf p$ be any extension of $\pi$ to a global generic.

Consider any definable subset $X$ of $B$ containing $B ^{00}$. Then finitely many right translates of $BaX$ cover $BaB$, so by $(*)$, $BaX$ is generic. Hence boundedly many left translates of $B^{00}aX$ cover $G$, which implies that one of them, say $gB^{00}aX$, is defined by a partial type which is contained in $\bf p$.

Since boundedly many right translates of $gB^{00}aB^{00}$ by elements from $X$ cover $gB^{00}aX$, one of them, say $gB^{00}aB^{00}x$ where $x \in X$, is defined by a partial type which is contained in $\bf p$.

By the definition of $\pi$, we have that $\bf p$ extends the partial type $x \in G^{00}$, so $gB^{00}aB^{00}x \cap G^{00} \ne \emptyset$. Hence $gx(x^{-1}B^{00}aB^{00}x) \cap G^{00} \ne \emptyset$. Thus $gx \in G^{00}$, which implies that $g=hx^{-1}$ for some $h \in G^{00}$.

So a partial type defining $B^{00}aX$ is contained in $g^{-1}{\bf p}=xh^{-1}{\bf p}=x{\bf p}$ (as $Stab({\bf p})=G^{00}$). On the other hand, a partial type defining $xB^{00}a^{-1}B^{00}$ is also contained in $x{\bf p}$. Hence 

$$(!)\;\;\;\;\;\; B^{00}aX \cap xB^{00}a^{-1}B^{00} \ne \emptyset.$$

So we have proved that for every definable subset $X$ of $B$ containing $B^{00}$ there is $x \in X$ such that $(!)$ holds. Thus by the compactness theorem,  $B^{00}aB^{00} \cap B^{00}a^{-1}B^{00} \ne \emptyset$ and hence $B^{00}aB^{00}=B^{00}a^{-1}B^{00}$. \hfill $\square$\\[2mm]
By the subclaim, there is $b \in B^{00}$ such that $(ab)^2 \in B$. Put $i:=ab$. Then, of course, $i \in G^{00}$. We will show that $i$ is an involution. 

Since $a$ is generic over a name of $B$ and $\uth(B)=1$, we get that $i \ne e$. Suppose for a contradiction that $i^2 \ne e$. As $i^2 \in B$, we have $B(i^2)=B$. But $i$ centralizes $i^2$, so $B^{i}=B(i^{2})^{i}=B(i^2)=B$. Thus $i \in N(B)$, which implies $a \in N(B)$, a contradiction with the fact that $a$ is generic over a name of $B$. \hfill $\square$\\

Notice that by Corollary \ref{involution corollary 2}, we get that all involutions have infinite centralizers. Indeed, if there was an involution $i$ with a finite centralizer, then $\uth(i^G)=2$, so there would be an involution of $\uth$-rank 2. Then by Corollary \ref{involution corollary 2}, we would get a non-trivial element whose centralizer has a finite index in $G$, a contradiction.

\begin{claim}
If $i$ and $j$ are involutions such that $B(i) \ne B(j)$, then $C(ij)$ is finite.
\end{claim}
\noindent {\em Proof.} The argument is similar as in the Morley rank 2 case. Indeed, we have $(ij)^i=(ij)^j = ji= (ij)^{-1}$. Hence $B(ij)^i=B((ij)^{-1})=B(ij)$, which implies that $i \in N(B(ij))$. Similarly, $j \in N(B(ij))$. But $i$ and $j$ have infinite centralizers. Hence by Claim 5, $i,j \in C(B(ij))$. Similarly, by Claim 5, we easily conclude that $i \in C(B(i))$ and $j \in C(B(j))$. In virtue of Claim 2(i), we get $B(i)=B(ij)=B(j)$, a contradiction. \hfill $\square$ 

\begin{claim} Every Borel $B$ contains at most finitely many involutions.
\end{claim}
\noindent {\em Proof.} Let $I$ be the set of all involutions. Suppose for a contradiction that $I \cap B$ is infinite. Then $I \cap B^g$ is also infinite for every $g \in G$. Now add a name of $B$ to the language. 

Take $g$ so that $N(B)g \notin \acl^{eq}(\emptyset)$. Let $\phi(x,y)$ be a formula over $\emptyset$ such that $\phi(\C,N(B)g)=I \cap B^{N(B)g}=I \cap B^g$.  
We see that for any $N(B)g_1 \ne N(B)g_2$ realizing $\tp(N(B)g)$ we have $\phi(\C,N(B)g_1) \cap \phi(\C,N(B)g_2)=I \cap B^{g_1} \cap B^{g_2}=\{ e \}$. 

So by Lemma \ref{finite intersections good enough}, there is an involution $b \in \phi(\C,N(B)g)=I \cap B^{g}$ such that $\uth(b)=2$, a contradiction with Corollary \ref{involution corollary 2}. \hfill $\square$\\

Let $X$ be the set of all elements of $G^{00}$ with finite centralizers. By Claim 4 and Claim 6, we see that $G^{00} \setminus X$ is a union of connected components of Borels.

\begin{claim}
(i) For every $a \in X$, if $a^n \in G^{00} \setminus X$, then $a^n =e$.\\
(ii) There is no $a \in X$ such that $a^2 \in G^{00} \setminus X$.\\
(iii) Every element of $X$ has an odd exponent.
\end{claim}
\noindent {\em Proof.} (i) Suppose $a \in X$, $a^n \in G^{00} \setminus X$ and $a^n \ne e$. Since $a^{-1}a^na = a$, we get $B(a^n)^a=B(a^n)$, so by Claim 6, we get $a \in N_{G^{00}}(B(a^n)^{00})=B(a^n)^{00}$, a contradiction.\\
(ii) If $a \in X$ and $a^2 \in G^{00} \setminus X$, then by (i), we see that $a^2=e$. So we get an involution with a finite centralizer, a contradiction.\\
(iii) Suppose that some $a \in X$ has an even exponent, say $2n$. Then $a^n \ne e$, so by (i), $a^n \in X$. But we also have $(a^n)^2 =e$, a contradiction with (ii). \hfill $\square$

\begin{claim}
For every $a \in X$ there is a unique $b \in G^{00}$ such that $b^2=a$. Moreover, $b \in \langle a \rangle$.
\end{claim}
\noindent {\em Proof.} By Claim 10(iii), we have that $\langle a \rangle$ is an abelian group of odd order, say $n$. So there is a unique $b \in \langle a \rangle$ such that $b^2=a$. 

Now suppose $b_1^2=a$ for some $b_1 \in G^{00}$. Then $b_1 \in X$. Since $b_1^{2n}=a^n=e$, we get $(b_1^n)^2=e$, so by Claim 10(i, ii), $b_1^n=e$. But $\langle a \rangle$ is a subgroup of $\langle b_1 \rangle$ of cardinality $n$. Thus $\langle a \rangle = \langle b_1 \rangle$, which implies $b=b_1$. \hfill $\square$\\

Now we will introduce and use a function similar to that appearing in \cite{involutions}, which comes from the theory of so-called black box groups.

By Claim 7, we can choose an involution $i \in G^{00}$. Let $B:=B(i)$. By Claim 6, we have $N_{G^{00}}(B^{00})=B^{00}$.

Notice that if $g \in G^{00} \setminus B$, then $B(i^g) \ne B(i)$, so by Claim 8, we get $ii^g \in X$. Hence $\sqrt{ii^g}$ is well-defined in $G^{00}$ by Claim 11. So we can define a function $f : G^{00} \setminus B \to G^{00}$ putting 

$$f(g)=\sqrt{ii^g}g^{-1}.$$ 
\begin{claim} $rng(f)=B^{00}$.
\end{claim}
\noindent {\em Proof.} First we check that $rng(f) \subseteq B^{00}$. By Claim 6, we get $B^{00}=C(i) \cap G^{00}$, so it is enough to show that for $g \in G^{00} \setminus B$, $i\sqrt{ii^g}g^{-1}=\sqrt{ii^g}g^{-1}i$. We have the following sequence of equivalent conditions: $i\sqrt{ii^g}g^{-1}=\sqrt{ii^g}g^{-1}i$ iff $i\sqrt{ii^g}i^g=\sqrt{ii^g}$ iff $i\sqrt{ii^g}i^gi\sqrt{ii^g}i^g=ii^g$ iff $\sqrt{ii^g}(ii^g)^{-1}\sqrt{ii^g}=e$. Since $ii^g$ and $\sqrt{ii^g}$ commute, we see that the last condition is true. So $rng(f) \subseteq B^{00}$.
Now we easily check that for $b \in B^{00}$ we have $f(bg)=f(g)b^{-1}$. So $rng(f)=B^{00}$. \hfill $\square$\\

Take any $g \in G^{00} \setminus B$. By Claim 12, we get $\sqrt{ii^g} =b_1g$ for some $b_1 \in B^{00}$. Hence $ii^g=b_1gb_1g$, so $g^{-1}=ib_1gb_1i$. Let $b=ib_1$. Since $i,b_1 \in B^{00}$, we get $b \in B^{00}$ and $(bg)^2=e$. But $g \notin B^{00}$, so $bg$ is an involution.

Hence we have proved the following statement:

\begin{enumerate}
\item[$(*)$] for every $g \in G^{00}$, $B^{00}g \cap I \ne \emptyset$ where $I$ is the set of involutions in $G^{00}$.
\end{enumerate}

\noindent
In other words, $B^{00}I=G^{00}$.
 
Notice that everything that we have proved about $B$ (including $(*)$) holds for any Borel whose connected component contains an involution. 

Let $d$ be a name of $B$.
Replacing $B$ by $B^h$ for some $h \in G^{00}$, if necessary, we can assume that  $d \notin \acl^{eq}(\emptyset)$. Now we can easily choose an element $g \in G^{00} \setminus (X \cup I \cup \{ e \})$ so that $g \thind d$. Then $g \in G^{00} \setminus B^{00}$. 

\begin{claim}
(i) $\uth(g^G) = 1$\\
(ii) $g^B$ is infinite.
\end{claim}
\noindent {\em Proof.} Point (i) follows from Lascar inequalities and the fact that $C(g)$ is infinite. 
To see (ii), notice that if $g^B$ were finite, then $C(g) \cap B$ would be infinite, so $C(g)\cap B=B$. Then $g \in G^{00} \cap C(B)=B^{00}$, a contradiction.
\nolinebreak \hfill $\square$ \\ 
  
Now add $g$ to the language. Then $d \notin \acl^{eq}(\emptyset)$ and there is a formula $\phi(x,y)$ over $\emptyset$ such that $\phi(\C,d)=g^B$.

\begin{claim}
If $d_1 \ne d_2$ are any realizations of $\tp(d)$, then $\phi(\C,d_1) \cap \phi(\C,d_2)$ is finite.
\end{claim}
\noindent {\em Proof.} There are automorphisms $f_1$ and $f_2$ such that $f_1(d)=d_1$ and $f_2(d)=d_2$.  Let $B_1=f_1[B]$ and $B_2=f_2[B]$. Then $\phi(\C,d_1)=g^{B_1}$, $\phi(\C,d_2)=g^{B_2}$, and $B_1 \ne B_2$. So by $(*)$ (applied to $B_1$ and $B_2$), we get that there are $b_1 \in B_1^{00}$, $b_2 \in B_2^{00}$, and $j_1,j_2 \in I$ such that $g=b_1j_1=b_2j_2$. Since $g \notin I$, we have $b_1, b_2 \ne e$. 

Suppose for a contradiction that $g^{B_1} \cap g^{B_2}$ is infinite. Then $b_1j_1^{B_1} \cap b_2j_2^{B_2}$ is infinite. So there are infinite subsets $\{k_n: n \in \omega\}$ and $\{ l_n : n \in \omega \}$ of $I$ such that $c:=b_2^{-1}b_1=k_nl_n$ for all $n \in \omega$. Since $B_1 \ne B_2$ and $b_1,b_2 \ne e$, we get $c \ne e$. 


Notice that there is $n \in \omega$ such that $B(k_n) \ne B(l_n)$. Otherwise, $k_n,l_n \in B(k_n)^{00}$ and so $c \in B(k_n)^{00}$ for all $n \in \omega$. 
Hence $k_n, l_n \in B(k_{0})^{00}$ for all $n \in \omega$, a contradiction with Claim 9. 

So by Claim 8, we get that $C(c)$ is finite. On the other hand, $c^{k_n}=(k_nl_n)^{k_n}=(k_nl_n)^{-1}=c^{-1}$, so $k_n \in N(C(c))$ for all $n\in \omega$, which implies that $N(C(c))$ is infinite. This is a contradiction. \hfill $\square$\\

By Lemma \ref{finite intersections good enough}, Claim 13(ii), and Claim 14, we get that there is an element $h \in \phi(\C,d)$ such that $\uth(h)=2$. But $h \in g^G$ which is $\emptyset$-definable (as we have added $g$ to the language). Thus $\uth(g^G)=2$, a contradiction to Claim 13(i). This completes the proof of Theorem \ref{3.1}. \hfill $\blacksquare$

\begin{remark}
    We may see, by examining the proof above, that not only is there a solvable group of finite index, but that we may take this group to be definable.  However, in any case, given a solvable group of finite index, we may find a definable solvable group also of finite index by the argument of, for instance, Theorem 3.17 of \cite{groupes stables}.
\end{remark}

Notice that Theorem 3.1 is true for an arbitrary independence relation $\starind$, because $\uth$-rank is less or equal to $U^*$-rank. 

From the proof of Theorem 3.1, we get that, after possibly passing to a definable subgroup of finite index and quotienting by its finite center, the group is solvable of solvability degree at most 2.

Let us finish this section with the following conjecture which generalizes Theorem  \ref{Main Theorem 2}.

\begin{conjecture} 
Suppose $G$ is superrosy and has NIP. Then if $G$, and every definable subgroup of $G$, is definably amenable (in the sense of \cite{NIP paper}) and of $\uth$-rank 2, it is solvable-by-finite.
\end{conjecture}

\section{Superrosy fields}

First we make several remarks about (absolute) connected components and generic types in fields. Then we adapt Macintyre's proof of \cite[Theorem 3.1]{groupes stables} to prove Theorem \ref{Main Theorem 3}. The last part of this section is devoted to the proof of Theorem \ref{Main Theorem 4}.

Suppose $K$ is an infinite field definable in a monster model $\C \models T$. The additive and multiplicative groups of $K$ will be denoted by $K^+$ and $K^*$, respectively.

\begin{proposition}
(i) If $(K^{+})^{00}$ exists (e.g. if $K^{+}$ has fsg or $T$ satisfies NIP), then $(K^{+})^{00}=K^{+}$.\\ 
(ii) If $K^{+}$ has fsg, then it has a unique global generic type $\bf p$ and $Stab_{K^*}({\bf p})=K^*=(K^{*})^{00}$.\\ 
(iii) If there is a global type $\bf p$ such that $Stab_{K^*}({\bf p})=K^*$ and $K^{*}$ has a global generic type (e.g. if $K^{*}$ has fsg), then $\bf p$ is the unique global generic type of $K^{*}$.
\end{proposition}
\noindent {\em Proof.}
(i) Since $(K^{+})^{00}$ is the smallest subgroup of bounded index in $K^+$, we see that for any $k \in K^{*}$, $k(K^{+})^{00}$ also has this property. Hence $k(K^+)^{00}=(K^+)^{00}$, which means that $(K^+)^{00}$ is a nontrivial ideal of $K$. Thus $(K^+)^{00}=K^{+}$.\\ 
(ii) Assume that $K^+$ has fsg. By (i) and Proposition \ref{unique generic}, we get the  existence and uniqueness of a global generic type ${\bf p}$ of $K^+$. Then we see that for every $k \in K^*$, $k{\bf p}$ is also generic, so $k{\bf p}={\bf p}$. Hence $Stab_{K^*}({\bf p})=K^*$.   On the other hand, it is clear that if $H$ is a type-definable subgroup of bounded index in $K^*$, then $Stab_{K^*}({\bf p}) \leq H$. So $(K^*)^{00}$ exists and is equal to $K^*$.\\
(iii) As in the proof of Proposition \ref{unique generic}, we show that every generic formula in $K^*$ belongs to ${\bf p}$. So, if there is a global generic type for $K^*$, it must be $\bf p$. \hfill $\blacksquare$   

\begin{proposition} Assume that $T$ is superrosy, $(K^{+})^{00}=K^{+}$ and
$(K^{*})^{00}=K^{*}$. Then for every $n >0$, the function $f(x)=x^n$ is onto and, if $char(K)=p$ is finite, the function $g(x)=x^p-x$ is also onto.
\end{proposition}
\noindent {\em Proof.} The proof is completely standard. Let us show the first part (the second one is similar). Let $a$ be $\thorn$-generic over $\emptyset$. Since $a \in \acl(a^n)$, Lascar inequalities give us that $a^n$ is also $\thorn$-generic. So we get that
$K^n$ has finite index in $K^{*}$, and by our assumption, we conclude that $K^n=K$. \hfill $\blacksquare$

\begin{proposition}
If $K^+$ has fsg and $L$ is a finite extension of $K$, then $L$ is definable in $\C$ and $L^{+}$ also has fsg. In particular, $(L^{+})^{00}=L^{+}$ and $(L^{*})^{00}=L^{*}$.  If $T$ is additionally superrosy, then $L$ satisfies the conclusion of Proposition 4.2.
\end{proposition} 
\noindent {\em Proof.}
Of course $L$ is definable in $K^{\times n}$ for some $n$. That $L^+$ has fsg follows by induction from \cite[Proposition 4.5]{NIP paper} and the fact that $L^{+}$ can be identified with $(K^{+})^{\times n}$. The rest is a consequence of Propositions 4.1 and 4.2. \hfill $\blacksquare$

\begin{proven main theorem}
Suppose that $T$ is superrosy and $K^+$ has fsg. Then $K$ is algebraically closed.
\end{proven main theorem}
\noindent
{\em Proof.} If not, then (as in Macintyre's proof, see \cite[Theorem 3.1]{groupes stables}) by Galois theory, there are finite extensions $K \subseteq L \subseteq F$ such that $F$ is a cyclic extension of $L$ of prime degree,  $q$, and if $q$ is different from the characteristic of $K$, then a primitive $q$th root of 1 belongs to $L$; then we get a contradiction with the last part of Proposition 4.3. \hfill $\blacksquare$\\

The following conjecture seems more difficult to prove. But it may be easier than the corresponding conjecture concerning supersimple fields, namely that each supersimple field is pseudo-algebraically closed.

\begin{conjecture}
Suppose $T$ is superrosy and has NIP. Then $K$ is either an algebraically closed or a real closed field.
\end{conjecture}

Suppose $G$ is a group definable in $\C \models T$. We will prove the following theorem which, by Corollary \ref{icc on centralizers}, implies Theorem \ref{Main Theorem 4}.

\begin{theorem} Assume that $G$ has hereditarily fsg,  the definable quotients of definable subgroups of $G$ satisfy icc on centralizers, $\uth(G)$=2 and $G$ is not nilpotent-by-finite. Then, after possibly passing to a definable subgroup of finite index and quotienting by its finite center, $G$ is (definably) the semidirect product of the additive and multiplicative groups of an algebraically closed field $F$  interpretable in $\C$, and moreover $G = G^{00}$.
\end{theorem}
\noindent
\noindent {\em Proof.} 
By Theorem \ref{3.1}, we know that $G$ is solvable-by-finite.  In fact, it has a definable solvable subgroup of finite index. 
We may assume that $G$ is centralizer connected, centerless and solvable. Let $U$ be a definable normal commutative subgroup of $G$ of $\uth$-rank $1$. We may assume that $U$ is centralizer connected in the sense of $G$, namely that for no $g\in G$ is $C(g)\cap U$ a proper subgroup of $U$ of finite index. Note that $C(U)$ 
(centralizer of $U$ in $G$) has infinite index in $G$ (otherwise, $G$ is nilpotent-by-finite). Also $G/C(U)$ being of $\uth$-rank $1$ is abelian-by-finite. It follows that there is $b\in G\setminus C(U)$ such that $C(b)$ is infinite (and so of $\uth$-rank $1$). (Otherwise, every conjugacy class in $G\setminus C(U)$ has $\uth$-rank $2$ so there are only finitely many of them, but then $G/C(U)$ has only finitely many conjugacy classes, contradicting it being infinite and abelian-by-finite.) By choice of $U$ and $b$, $C(b)\cap U$ is finite, hence $C(b)U$ (the group generated by $C(b)$ and $U$ which is definable) has $\uth$-rank $2$ so finite index in $G$. We claim that $C(b)\cap U = \{e\}$. For otherwise, 
the centralizer of every element in $C(b) \cap U \ne \{ e\}$ is of $\uth$-rank 2, which contradicts $G$ being centralizer connected and centerless.
Let $T$ be a commutative definable subgroup of $C(b)$ of finite index. It follows
likewise that $T\cap C(U)$ is also trivial and for every $u \in U \setminus \{e\}$, $T \cap C(u)$ is finite. Now the group $UT$ has $\uth$-rank $2$, so we may assume it equals $G$. So to summarize the situation we have:

\begin{enumerate}
\item[$(*)$] $G = UT$ (semidirect product of $U$ and $T$), $U, T$ are commutative of $\uth$-rank $1$, $C(U)\cap T = \{e\}$ and for every $u \in U \setminus \{e\}$, $T \cap C(u)$ is finite. 
\end{enumerate}

We will write $U$ additively, and $T$ multiplicatively. $T$ acts by conjugation on $U$, and we sometimes let $t\cdot u$ denote $u^{t}=tut^{-1}$. By the last part of $(*)$, each orbit except $\{ 0 \}$ is infinite. So by $(*)$, $T$ acts regularly on each orbit except $\{ 0 \}$. (Otherwise, there are $u \in U \setminus \{ 0 \}$ and $t \in T \setminus \{ e \}$ such that $u^t=u$. Then every element of the infinite orbit $u^T$ is stabilized by $t$, so $C(t) \cap U$ is infinite, and hence $t \in T \cap C(U)$, a contradiction.) Hence as $U$ has $\uth$-rank $1$, there are only finitely many orbits. But $U^{00}$ is clearly a union of such orbits. Hence $U^{00}$ is definable. So we may assume:
\begin{enumerate}
\item[$(**)$] $U = U^{00}$ (hence we know there is a unique generic type in $U$). 
\end{enumerate}
Note that all our previous assumptions remain valid. In particular $T$ acts freely on $U\setminus\{0\}$. Let $C_{1},\dots,C_{r}$ be the orbits of $U\setminus\{0\}$ under $T$. By $(**)$, exactly one of them, without loss $C_{1}$ is generic in $U$. If it so happened that $r=1$ then we easily get our desired conclusion (as in the 
finite Morley rank case). In fact we will undertake an analysis which will eventually show that $r=1$ anyway.

The important thing we will use is that each element of $U$ is a sum of generics and thus a sum of elements of $C_{1}$. Each element of $T$ defines an endomorphism (in fact automorphism) of $U$. Let $S$ be the ring of endomorphisms of $U$ generated by $T$. Note that two different elements of $T$ define different endomorphisms (in fact automorphisms) of $U$, so we can and will identify an element of $T$ with the endomorphism it defines. Moreover, multiplication in $T$ is just the restriction to $T$ of multiplication (composition) in $S$. As $T$ is commutative, so is $S$. For $s\in S$ we still write the action on $U$ as $\cdot$. We write $T+T$ for the subset of $S$ consisting of elements $s+t$ for $s,t\in T$.
Fix an element $u\in C_{1}$.\\[2mm]
{\bf Claim} {\em (i) Any $s\in S$ is determined by $s\cdot u$.
\newline
(ii) $S = T+T$.
\newline 
(iii) The ring $S$ is an (interpretable) field.
\newline
(iv) Let $i:S\to U$ be given by $i(s) = s\cdot u$. Then $i$ induces an isomorphism between the field $S$ and $(U,+,\otimes)$ (some definable $\otimes$) and moreover this field is algebraically closed.}\\[2mm]
\noindent {\em Proof.} This is actually quite routine and implicit or explicit in the literature, but we will give some details anyway.
\newline
(i) As every element of $U$ is a sum of elements of $C_{1}$ and as $C_{1} = T\cdot u$, we see that every element of $U$ is of the form $s\cdot u$ for some $s\in T+T$. Let $s_{1},s_{2}\in S$, and suppose that $s_{1}\cdot u = s_{2}\cdot u$. Let $x\in U$, and suppose $s\in T+T$ is such that $x = s\cdot u$. Then
$s_{1}\cdot x = s_{1}\cdot(s\cdot u) = s\cdot(s_{1}\cdot u)) = s\cdot(s_{2}\cdot u) = s_{2}\cdot(s\cdot u) = s_{2}\cdot x$. 
\newline
(ii) follows from (i) as $(T+T)\cdot u = U$. 
\newline 
(iii) Suppose first that $s\in S$, $x\in U$ and $s\cdot x = 0$. Now $x\in C_{i}$ for some $i=1,\dots,r$ and $C_{i} = T\cdot x$. As $S$ is commutative, $s$ is $0$ on $C_{i}$, so $ker(s)$ (a definable subgroup of $U$) is infinite, thus as $U$ is connected of $\uth$-rank $1$, $ker(s) = U$ and $s = 0$. So we have shown that any 
nonzero $s\in S$ is injective, hence also surjective (by connectedness of $U$ again). The existence of inverses follows easily: let $s\in S$ be nonzero. From what we have seen let $x\in U$ be such that $s\cdot x = u$. Let $s'\in S$ be such that $s'\cdot u = x$. So $s\cdot (s'\cdot u) = u$, and by part (i), $s\cdot s'$ is the identity. So $s'$ is the inverse of $s$. We have shown that $S$ is a field. 
Identifying $S$ with $T\times T/E$ for a suitable definable equivalence relation $E$ (using (i) and (ii)), we see that addition and multiplication are definable.

So we have an interpretable definable field. But we only have a theorem telling us the structure of such a field when the underlying additive group has fsg. This is the point of:
\newline
(iv) We know by (i) that $i$ is a bijection. Moreover $i$ clearly takes addition on $S$ to addition on $U$. So we have a definable field structure $(U,+,\otimes)$ on $U$ whose additive part has the fsg (by our hypothesis on $G$). By Theorem \ref{Main Theorem 3} the field is algebraically closed. \hfill $\square$\\  
%

Let us now complete the proof of the theorem. Consider the definable field structure $F = (U,+,\otimes)$ expanding $(U,+)$. From Proposition 4.1 $(U^{*},\otimes)$ is (absolutely) connected, in particular has no proper definable 
subgroup of finite index. But $T$ embeds definably in $(U^{*},\otimes)$ via $t\to t\cdot u$, and $T$ as well as $U^{*}$ has $\uth$-rank $1$. This forces $T\cdot u$ to equal $U\setminus \{0\}$. So in fact the action of $T$ on $U$ is the action of $F^{*}$ on $F^{+}$. 

For absolute connectedness, we have that $U = U^{00}$ and that $T 
= T^{00}$ (being isomorphic to the multiplicative group of $F$). It follows that $G$, the semidirect product of $U$ and $T$ is absolutely connected: $G^{00} \cap U=U$ and $G^{00} \cap T=T$, hence $G=G^{00}$. \hfill $\blacksquare$

\end{document}